\documentclass[12pt]{article}

\usepackage{amsfonts,amssymb,subeqnarray,eqsection,indent}
\usepackage{e-jc}

\date{\dateline{September 28, 2009}{January 29, 2010}
      \small  05C15, 05C30, 05C35, 05C40, 82B20, 90B10 (Secondary)}

\begin{document}

\title{Maxmaxflow and Counting Subgraphs}

\author{
  {\small Bill Jackson}                                    \\[-2mm]
  {\small\it School of Mathematical Sciences}  \\[-2mm]
  {\small\it Queen Mary University of London} \\[-2mm]
  {\small\it Mile End Road} \\[-2mm]
  {\small\it London E1 4NS, England}                         \\[-2mm]
  {\small\tt B.JACKSON@QMUL.AC.UK}                        \\[5mm]
  {\small Alan D.~Sokal\thanks{Also at Department of Mathematics,
           University College London, London WC1E 6BT, England.}}   \\[-2mm]
  {\small\it Department of Physics}       \\[-2mm]
  {\small\it New York University}         \\[-2mm]
  {\small\it 4 Washington Place}          \\[-2mm]
  {\small\it New York, NY 10003 USA}      \\[-2mm]
  {\small\tt SOKAL@NYU.EDU}               \\[-2mm]
  {\protect\makebox[5in]{\quad}}  
  \\
}
\maketitle
\thispagestyle{empty}   

\begin{abstract}
We introduce a new graph invariant $\Lambda(G)$ that we call maxmaxflow,
and put it in the context of some other well-known graph invariants,
notably maximum degree and its relatives.
We prove the equivalence of two ``dual'' definitions of maxmaxflow:
one in terms of flows, the other in terms of cocycle bases.
We then show how to bound the total number
(or more generally, total weight)
of various classes of subgraphs of $G$
in terms of either maximum degree or maxmaxflow.
Our results are motivated by a conjecture
that the modulus of the roots of the chromatic polynomial
of $G$ can be bounded above by a function of $\Lambda(G)$.
\end{abstract}



\clearpage

\newtheorem{defin}{Definition}[section]
\newtheorem{definition}{Definition}[section]
\newtheorem{prop}[defin]{Proposition}
\newtheorem{proposition}[defin]{Proposition}
\newtheorem{lem}[defin]{Lemma}
\newtheorem{lemma}[defin]{Lemma}
\newtheorem{guess}[defin]{Conjecture}
\newtheorem{ques}[defin]{Question}
\newtheorem{question}[defin]{Question}
\newtheorem{prob}[defin]{Problem}
\newtheorem{problem}[defin]{Problem}
\newtheorem{thm}[defin]{Theorem}
\newtheorem{theorem}[defin]{Theorem}
\newtheorem{cor}[defin]{Corollary}
\newtheorem{corollary}[defin]{Corollary}
\newtheorem{conj}[defin]{Conjecture}
\newtheorem{conjecture}[defin]{Conjecture}

\newtheorem{pro}{Problem}
\newtheorem{clm}{Claim}
\newtheorem{con}{Conjecture}

%
%
\newcounter{example}[section]
\newenvironment{example}%
{\refstepcounter{example}
 \bigskip\par\noindent{\bf Example \thesection.\arabic{example}.}\quad
}%
{\quad $\Box$}
\def\bexam{\begin{example}}
\def\eexam{\end{example}}

\renewcommand{\theenumi}{\alph{enumi}}
\renewcommand{\labelenumi}{(\theenumi)}
\def\prf{\par\noindent{\bf Proof.\enspace}\rm}
\def\rmk{\par\medskip\noindent{\bf Remark.\enspace}\rm}

\newcommand{\be}{\begin{equation}}
\newcommand{\ee}{\end{equation}}
\newcommand{\<}{\langle}
\renewcommand{\>}{\rangle}
\newcommand{\widebar}{\overline}
\def\reff#1{(\protect\ref{#1})}
\def\spose#1{\hbox to 0pt{#1\hss}}
\def\ltapprox{\mathrel{\spose{\lower 3pt\hbox{$\mathchar"218$}}
 \raise 2.0pt\hbox{$\mathchar"13C$}}}
\def\gtapprox{\mathrel{\spose{\lower 3pt\hbox{$\mathchar"218$}}
 \raise 2.0pt\hbox{$\mathchar"13E$}}}
\def\textprime{${}^\prime$}
\def\proof{\par\medskip\noindent{\sc Proof.\ }}
\newcommand{\qed}{\quad $\Box$ \medskip \medskip}
\def\proofof#1{\bigskip\noindent{\sc Proof of #1.\ }}
\def\half{ {1 \over 2} }
\def\third{ {1 \over 3} }
\def\twothird{ {2 \over 3} }
\def\smfrac#1#2{\textstyle{#1\over #2}}
\def\smhalf{ \smfrac{1}{2} }
\newcommand{\real}{\mathop{\rm Re}\nolimits}
\renewcommand{\Re}{\mathop{\rm Re}\nolimits}
\newcommand{\imag}{\mathop{\rm Im}\nolimits}
\renewcommand{\Im}{\mathop{\rm Im}\nolimits}
\newcommand{\sgn}{\mathop{\rm sgn}\nolimits}
\def\hboxscript#1{ {\hbox{\scriptsize\em #1}} }

\newcommand{\restrict}{\upharpoonright}
\renewcommand{\emptyset}{\varnothing}

\newcommand{\scra}{{\mathcal{A}}}
\newcommand{\scrb}{{\mathcal{B}}}
\newcommand{\scrc}{{\mathcal{C}}}
\newcommand{\scrf}{{\mathcal{F}}}
\newcommand{\scrg}{{\mathcal{G}}}
\newcommand{\scrh}{{\mathcal{H}}}
\newcommand{\scrl}{{\mathcal{L}}}
\newcommand{\scro}{{\mathcal{O}}}
\newcommand{\scrp}{{\mathcal{P}}}
\newcommand{\scrr}{{\mathcal{R}}}
\newcommand{\scrs}{{\mathcal{S}}}
\newcommand{\scrt}{{\mathcal{T}}}
\newcommand{\scrv}{{\mathcal{V}}}
\newcommand{\scrw}{{\mathcal{W}}}
\newcommand{\scrz}{{\mathcal{Z}}}
\newcommand{\scrbt}{{\mathcal{BT}}}
\newcommand{\scrbf}{{\mathcal{BF}}}


\newenvironment{sarray}{
      \textfont0=\scriptfont0
      \scriptfont0=\scriptscriptfont0
      \textfont1=\scriptfont1
      \scriptfont1=\scriptscriptfont1
      \textfont2=\scriptfont2
      \scriptfont2=\scriptscriptfont2
      \textfont3=\scriptfont3
      \scriptfont3=\scriptscriptfont3
    \renewcommand{\arraystretch}{0.7}
    \begin{array}{l}}{\end{array}}

\newenvironment{scarray}{
      \textfont0=\scriptfont0
      \scriptfont0=\scriptscriptfont0
      \textfont1=\scriptfont1
      \scriptfont1=\scriptscriptfont1
      \textfont2=\scriptfont2
      \scriptfont2=\scriptscriptfont2
      \textfont3=\scriptfont3
      \scriptfont3=\scriptscriptfont3
    \renewcommand{\arraystretch}{0.7}
    \begin{array}{c}}{\end{array}}

\tableofcontents

\clearpage

\section{Introduction}

An elementary result on graph colouring is that
the chromatic number $\chi(G)$ of a graph $G$
is at most one more than the maximum degree $\Delta(G)$.
A much deeper result is that the modulus of the roots (real or complex)
of the chromatic polynomial of $G$ can be bounded above by a
linear function of $\Delta(G)$, see \cite{Sokal_chromatic_bounds}.
Indeed, a similar bound holds when the maximum degree $\Delta(G)$
is replaced by the second-largest degree $\Delta_2(G)$,
although the currently available proof of this fact
\cite[Corollary~6.4]{Sokal_chromatic_bounds}
is somewhat {\em ad hoc}\/.\footnote{
   Note that it is {\em not}\/ possible to go farther and obtain a bound
   in terms of the third-largest degree $\Delta_3$,
   as the chromatic roots of the generalized theta graphs $\Theta^{(s,p)}$
   --- which have $\Delta = \Delta_2 = p$ but $\Delta_3 = 2$ ---
   are dense in the whole complex plane with the possible exception
   of the disc $|q-1| < 1$
   \cite[Theorems 1.1--1.4]{Sokal_chromatic_roots}.
}

One obvious drawback in all these results is that we can make
the maximum degree and second-largest degree arbitrarily large
by gluing together many copies of $G$ in a tree-like fashion
at cut vertices, without changing the chromatic number or
the chromatic roots.  Another (related) drawback is that there is no
obvious way to extend these results from graphs to matroids
and thereby to obtain dual results for nowhere-zero flows and
the roots of flow polynomials.

The purpose of this paper is to introduce a new graph invariant
$\Lambda(G)$ that we call {\em maxmaxflow}\/, which we conjecture
will give a more natural upper bound on chromatic roots.
The maxmaxflow $\Lambda(G)$ is defined as the maximum,
over all pairs of distinct vertices $x,y$ of $G$,
of the maximum number of pairwise edge-disjoint $xy$-paths.
It is easy to see that $\Lambda(G)$ is less than or equal to $\Delta_2(G)$,
and that the maxmaxflow of any graph is equal to the largest maxmaxflow
in its blocks (maximal non-separable subgraphs).
We will show that $\Lambda(G)$ can equivalently be defined
in terms of the bases of the cocycle space of $G$,
so that the definition of maxmaxflow can be extended to binary matroids.
We will furthermore see that $\Lambda(G)$ is at
least as large as the degeneracy number $D(G)$ of $G$,
so that we have $\chi(G)\leq D(G)+1 \leq \Lambda(G)+1$.
We conjecture that $\Lambda(G)$ can also be used to give a bound
on the chromatic roots of $G$:

\begin{conjecture}
   \label{conj.chrom.1}
There exist universal constants $C(\Lambda) < \infty$ such that all
the chromatic roots (real or complex) of all loopless graphs of
maxmaxflow $\Lambda$ lie in the disc $|q| \le C(\Lambda)$. Indeed,
we conjecture that $C(\Lambda)$ can be taken to be linear in
$\Lambda$.
\end{conjecture}

\noindent
This conjecture first appeared in \cite[Section~7]{Sokal_chromatic_bounds}
and was inspired by a suggestion of Shrock and Tsai
\cite{Shrock_98e,Shrock_99a}.
It has very recently been proven for series-parallel graphs
by Royle and Sokal \cite{Royle-Sokal}.

An important step in the proof \cite{Sokal_chromatic_bounds}
that the chromatic roots of $G$ can be bounded in terms of $\Delta(G)$
is obtaining an exponential upper bound in terms of $\Delta(G)$
for the number of connected $m$-edge subgraphs
containing a fixed vertex of $G$.
The approach in \cite{Sokal_chromatic_bounds} is to
decompose a spanning subgraph of $G$ into its connected components
and to treat these components as a ``polymer gas''.   The desired
bound on chromatic roots then follows from standard bounds on the
zeros of a polymer-gas partition function, once one has the
exponential bound on the number of
connected $m$-edge subgraphs containing a specified vertex.

Unfortunately, the number of connected $m$-edge subgraphs
containing a fixed vertex {\em cannot}\/ be bounded in terms of $\Lambda(G)$.
This can easily be seen by taking $G$ to be large star:
we have $\Lambda(G)=1$ and yet there is no bound on the number of
connected $m$-edge subgraphs containing the central vertex.

Since both the chromatic polynomial and maxmaxflow ``factorize over blocks'',
it is natural to try to prove Conjecture~\ref{conj.chrom.1}
by modifying the arguments of \cite{Sokal_chromatic_bounds}
to decompose a spanning subgraph of $G$ into its blocks
rather than its connected components.
The main result of this paper, Corollary~\ref{cor.block.1},
is a first step in this direction.
It shows --- a result that some readers may find surprising ---
that the number of non-separable $m$-edge subgraphs
containing a fixed edge of $G$
satisfies an exponential upper bound in terms of $\Lambda(G)$.
This will be good enough to prove
Conjecture~\ref{conj.chrom.1} provided that other difficulties
(such as controlling the interaction between blocks) can be overcome.

Irrespective of the potential application to bounding chromatic roots,
we think that maxmaxflow is a natural graph invariant
that deserves further study
and that bounds on the number of subgraphs of various kinds
in terms of $\Delta(G)$ or $\Lambda(G)$ are of independent interest.\footnote{
   See Section~\ref{sec.maxmax} below for references to
   scattered earlier work concerning maxmaxflow.
}

The plan of this paper is as follows:
In Section \ref{sec.maxmax} we introduce maxmaxflow
and put it in the context of some other well-known graph invariants
(notably maximum degree and its relatives and degeneracy number).
In Section~\ref{sec.cocycle} we analyze cocycle bases and prove
the equivalence of the two definitions of maxmaxflow;
an important role in this proof is played by
{\em Gomory--Hu trees}\/ \cite{GH}.
The remainder of the paper is devoted to bounding
the total number (or more generally, total weight)
of various classes of subgraphs
in terms of either maximum degree or maxmaxflow.
Our basic approach is to start with a bound
(sometimes a known one, sometimes a new one) in terms of maximum degree,
and then see whether we can find a similar bound in terms of maxmaxflow.
After some brief preliminaries (Section~\ref{sec.prelim}),
we analyze walks and paths (Section~\ref{sec.walks})
and then trees and forests (Section~\ref{sec.trees}).
In Section~\ref{sec.conn} we consider connected subgraphs
and in Section~\ref{sec.block} we consider non-separable subgraphs.
Roughly speaking, the (more difficult) proofs in the later sections
are constructed by adapting ideas from the (easier) proofs
in the earlier sections.
We hope that, by organizing the paper
in terms of gradually increasing complexity of proof,
we have helped to reduce the mental burden on the reader.

\section{Maxmaxflow}
\label{sec.maxmax}

Let $G$ be a finite undirected graph with vertex set $V(G)$
and edge set $E(G)$;  in this paper all graphs are assumed
to be loopless, but multiple edges
are allowed unless explicitly specified otherwise. We shall say that
$G$ is {\em simple}\/ if it has no multiple edges.
Let $\Delta(G) = \max_{x \in V(G)} d_G(x)$ be the maximum degree of $G$,
and more generally let $\Delta_k(G)$ be the $k$th largest degree of $G$:
\be
   \Delta_k(G)   \;=\;
       \min\limits_{x_1,\ldots,x_{k-1} \in V(G)}
       \quad
       \max\limits_{x \in V(G) \backslash \{x_1,\ldots,x_{k-1}\}}
    d_G(x)
   \;.
\ee
We trivially have
\be
   \delta(G) \;\equiv\;
   \Delta_n(G) \;\le\; \cdots \;\le\; \Delta_3(G) \;\le\; \Delta_2(G)
   \;\le\; \Delta_1(G)  \;\equiv\;  \Delta(G)  \qquad
\ee
where $n = |V(G)|$.
A special role will be played in this paper
by the second-largest degree, $\Delta_2(G)$.

For $x,y \in V(G)$ with $x \neq y$, {\em the maximum flow} from $x$
to $y$ in $G$ is
\begin{subeqnarray}
   \lambda_G(x,y)
      & = &  \hbox{max \# of edge-disjoint paths from $x$ to $y$} \\
      & = &  \hbox{min \# of edges separating $x$ from $y$}
\end{subeqnarray}
We then define the {\em maxmaxflow} of $G$\/ \be
   \Lambda(G)   \;=\;
   \max\limits_{\begin{scarray}
           x,y \in V(G) \\
           x \neq y
        \end{scarray}}
   \lambda_G(x,y)   \;.
\ee
[Note the contrast with the edge-connectivity, which is the {\em minimum}\/
of $\lambda_G(x,y)$ over $x \neq y$.]
Clearly $\lambda_G(x,y) \le \min[d_G(x), d_G(y)]$,
so that
\be
   \Lambda(G) \;\le\; \Delta_2(G)  \;.
\ee
We will show later (Proposition~\ref{prop2.10})
that $\Lambda(G) \ge \Delta_{n-1}(G)$.
Note that several cases can arise:
\begin{itemize}
   \item[(a)]  $\Lambda(G) = \Delta_2(G) = \Delta(G)$.
      Indeed, in any regular graph one has
      $\Lambda(G) = \Delta_i(G)$ for all $i$ ($1 \le i \le n$).
   \item[(b)]  $\Lambda(G) = \Delta_2(G) \ll \Delta(G)$.
      This occurs, for example,
      in stars $K_{1,r}$ and wheels $K_1 + C_r$.
   \item[(c)]  More generally, one can have
      $\Lambda(G) = \Delta_{j+1}(G) \ll \Delta_{j}(G)$
      for any fixed integer $j$. Moreover, such examples can be taken
      to be $k$-connected for arbitrarily large $k$.\footnote{
    {\sc Proof.}  For $1\leq i\leq j$, let $H_i$ be a $k$-connected graph
    with one vertex $v_i$ of degree $\Delta \gg k$
    and all other vertices of degree $k$.
    [Such graphs can be constructed by taking
    a $(k-1)$-connected $(k-1)$-regular graph with a large number $\Delta$
    of vertices and adding a new vertex $v_i$ adjacent to every other vertex.]
    Construct $G$ from the disjoint union of $H_1,H_2,\ldots,H_j$
    by adding $k$ edges between
    each pair $H_i-v_i$ and $H_{i+1}-v_{i+1}$ ($1\leq i\leq j-1$)
    in such a way that the set of edges of $G$ which do not belong to
    any $H_i$ are independent.
    [This can be done as long as  $|V(H_i)|\geq 2k+1$.]
     Then $G$ is $k$-connected and satisfies
     $\Lambda(G) = \Delta_{j+1}(G) =k+1$.
    [Since all pairs of vertices of $G$ with degrees greater than $k+1$
     are of the form $v_s,v_t$ with $s \neq t$,
     and hence are separated by a set of $k$ edges,
     we have $\Lambda(G) \le k+1$.
     On the other hand, if we choose two vertices $x,y \in H_1$
     that are both adjacent to $H_2$, we can find $k$ edge-disjoint
     $xy$-paths in $H_1$ (since $H_1$ is $k$-connected)
     and an extra $xy$-path passing through $H_2$.]
    But $\Delta_{j}(G)= \Delta$.
}
\end{itemize}

Note also that maxmaxflow has a naturalness property that
maximum degree and $k$th-largest degree lack, namely,
it ``trivializes over blocks'':
$\Lambda(G) = \max_{1 \le i \le b} \Lambda(G_i)$
where $G_1,\ldots,G_b$ are the blocks of $G$
(Proposition~\ref{prop2.triv_blocks}).

Maxmaxflow appears to have been considered sporadically
in the graph-theoretic literature.
Bollob\'as \cite[section I.5]{Bollobas_78} addresses
some extremal problems involving maxmaxflow in simple graphs
(he uses the term ``maximum local edge-connectivity'' and denotes it
 $\bar{\lambda}(G)$);
see likewise Mader \cite[section IV]{Mader_79}.
In particular, Mader \cite{Mader_73} has shown that
whenever an $n$-vertex graph has more than $k(n-1)/2$ edges,
it has maxmaxflow at least $k$,
but that for every $n \ge k \ge 2$
there exists an $n$-vertex graph
with exactly $\lfloor k(n-1)/2 \rfloor$ edges
and maxmaxflow $k-1$.\footnote{
   For $k=2,3$ this is easy.
   For $k=4$ it was proven earlier by Bollob\'as \cite{Bollobas_66},
   and for $k=5,6$ by Leonard \cite{Leonard_72,Leonard_73}.
}


An apparently very different quantity can be defined via cocycle bases.
For $X,Y$ disjoint subsets of $V(G)$,
let $E(X,Y)$ denote the set of edges in $G$ between $X$ and $Y$.
A {\em cocycle}\/ of $G$ is a set $E(X,X^c)$
where $X \subseteq V(G)$ and $X^c \equiv V(G) \backslash X$.
It is well-known that the cocycles of $G$
form a vector space over GF(2) with respect to symmetric difference;
this is called the {\em cocycle space}\/ of $G$.
Let $\widetilde{\Lambda}(G)$ be the minmax cardinality of the
cocycles in a basis, i.e.
\be
   \widetilde{\Lambda}(G)  \;=\;
   \min\limits_{\scrb}  \max\limits_{C\in \scrb} |C|
\ee
where the min runs over all bases $\scrb$ of the cocycle space of $G$.
Since one special class of cocycle bases consists
of taking the stars
$C(x) = E(\{x\}, \{x\}^c)$
for all but one of the vertices in each component of $G$, we clearly have
\be
   \widetilde{\Lambda}(G) \;\le\; \Delta_2(G)  \;.
\ee

The relationship, if any, between maxmaxflow and cocycle bases
is perhaps not obvious at first sight.
But we shall prove (Corollary~\ref{cordlambdatwiddle}) that
\be
   \Lambda(G)  \;=\;  \widetilde{\Lambda}(G) \;.
\ee
The two definitions thus give dual approaches to the same quantity.

Finally, define the {\em degeneracy number}\/
$D(G) = \max_{H \subseteq G} \delta(H)$,
where the max runs over all subgraphs $H$ of $G$,
and $\delta(H)$ is the minimum degree of $H$.
It is easy to see that
\be
   D(G)   \;\le\;   \Delta_2(G)
\ee
[if $H$ has at least two vertices, then $\delta(H) \le \Delta_2(G)$;
 otherwise $\delta(H)=0$].
%
%
We shall in fact show (Proposition~\ref{prop2.10}) that
\be
   D(G)   \;\le\;   \Lambda(G)   \;.
\ee
In summary, therefore, we have
\be
   D(G) \;\le\; \Lambda(G)  \;=\; \widetilde{\Lambda}(G)
    \;\le\; \Delta_2(G)  \;\le\; \Delta(G)
   \;.
\ee

The natural setting for the results of this paper is, in fact,
that of a finite undirected loopless (multi)graph $G$ equipped with
nonnegative real edge weights ${\bf w} = \{ w_e \} _{e \in E(G)}$.
Indeed, all of the aforementioned invariants have natural generalizations
to this context.  Define first the weighted degree of a vertex,
\be
   d_G(x,{\bf w})   \;=\;   \sum_{e \ni x} w_e
   \;.
\ee
We then set
\begin{eqnarray}
   \Delta(G,{\bf w})   & = &  \max\limits_{x \in V(G)} d_G(x,{\bf w})  \\[2mm]
   \Delta_k(G,{\bf w})   & = &
       \min\limits_{x_1,\ldots,x_{k-1} \in V(G)}
       \quad
       \max\limits_{x \in V(G) \backslash \{x_1,\ldots,x_{k-1}\}}
    d_G(x,{\bf w})
   \\[2mm]
   \delta(G,{\bf w})   & = &  \min\limits_{x \in V(G)} d_G(x,{\bf w})  \\[2mm]
   \widetilde{\Lambda}(G,{\bf w})  & = &
       \min\limits_{\scrb}  \max\limits_{C\in \scrb} \sum\limits_{e \in C} w_e
   \label{def_weighted_cocycle} \\[2mm]
   D(G,{\bf w})   & = &
       \max\limits_{H \subseteq G} \delta(H,{\bf w}|_H)
\end{eqnarray}
Likewise, max-flow quantities are naturally defined when the
$\{ w_e \}$ are interpreted as edge capacities:
\begin{subeqnarray}
   \lambda_G(x,y; {\bf w})
      & = &  \hbox{max flow from $x$ to $y$ with edge capacities ${\bf w}$} \\
      & = &  \hbox{min cut between $x$ and $y$ with edge capacities ${\bf w}$}
      \qquad
 \label{def_weighted_flow}
\end{subeqnarray}
and thence
\be
   \Lambda(G,{\bf w})   \;=\;
   \max\limits_{\begin{scarray}
           x,y \in V(G) \\
           x \neq y
        \end{scarray}}
   \lambda_G(x,y; {\bf w})   \;.
 \label{def_weighted_maxmaxflow}
\ee
In this generality we shall prove
\be\label{dlambdatwiddle}
  D(G,{\bf w}) \;\le\; \Lambda(G,{\bf w})  \;=\; \widetilde{\Lambda}(G,{\bf w})
    \;\le\; \Delta_2(G,{\bf w})  \;\le\; \Delta(G,{\bf w})
   \;.
\ee
The unweighted case corresponds to setting all edge weights to 1.

Let us make a remark about the treatment of multiple edges.
It is easy to see that all the quantities appearing in \reff{dlambdatwiddle}
are unchanged if we replace a family $e_1,\ldots,e_n$ of parallel edges
with weights $w_{e_1}, \ldots, w_{e_n}$ by a single edge $e$
with weight $w_e = \sum_{i=1}^n w_{e_i}$.
So, in proving \reff{dlambdatwiddle}, we could, if we wanted,
restrict attention to simple graphs;  but we don't bother,
because no simplification of the proof is obtained by doing so.
Likewise, the weighted counts discussed in
Sections~\ref{sec.walks} and \ref{sec.trees}
are unchanged by this replacement,
because the subgraphs in question (walks, paths, trees and forests)
can include at most one edge from a family of parallel edges.
So it would suffice to prove the bounds
in Sections~\ref{sec.walks} and \ref{sec.trees}
for simple graphs;  but once again, we refrain from making this assumption
because nothing is gained by doing so.
For the weighted counts discussed
in Sections~\ref{sec.conn} and \ref{sec.block}, by contrast,
no simple reduction of multiple edges can be performed,
because the subgraphs in question
{\em do}\/ permit the inclusion of multiple edges.
We shall therefore have to deal there with multigraphs in all our arguments.

\section{Cocycle Bases and Maxmaxflow} \label{sec.cocycle}

Given a graph $G$ and disjoint subsets $X,Y \subseteq V(G)$,
let $E(X,Y)$ denote the set of edges in $G$ between $X$ and $Y$.
A {\em cocycle}\/ of $G$ is a set $E(X,Y)$ where $X,Y$ is a
bipartition of $V(G)$;
note that $X = \emptyset$ and $Y = \emptyset$ are allowed.
Let $\oplus$ denote symmetric difference.
The following lemma is well known:

\begin{lem}\label{sum}
Let $C_1=E(X_1,Y_1)$ and $C_2=E(X_2,Y_2)$ be two cocycles in $G$. Then
$C_1\oplus C_2=E\left((X_1\cap X_2)\cup (Y_1\cap Y_2), \,
(X_1\cap Y_2)\cup (Y_1\cap X_2)\right)$.
\end{lem}

It follows that the set of all cocycles of $G$
forms a vector space over GF(2) with respect to symmetric difference.
This is the {\em cocycle space}\/ of $G$. Its dimension is
$|V(G)|-c(G)$, where $c(G)$ denotes the number of components of $G$.

\begin{lem}\label{unique}
Let $G$ be a connected graph and let $C$ be a cocycle of $G$.
Then $C$ corresponds to a unique bipartition of $V(G)$.
\end{lem}

\proof
Suppose $C=E(X_1,Y_1)= E(X_2,Y_2)$.
Since $C\oplus C=\emptyset$ there are no edges in $G$ from
$(X_1\cap X_2)\cup (Y_1\cap Y_2)$ to $ (X_1\cap Y_2)\cup (Y_1\cap X_2)$.
Since $G$ is connected, it follows that either
$(X_1\cap X_2)\cup (Y_1\cap Y_2)=\emptyset$
and hence $(X_1,Y_1)=(Y_2,X_2)$,
or else $ (X_1\cap Y_2)\cup (Y_1\cap X_2)=\emptyset$
and hence $(X_1,Y_1)=(X_2,Y_2)$.
\qed

\begin{lem}\label{sameside}
Let $G$ be a connected graph, let $C_1,C_2$ be  cocycles of $G$,
and let $x,y$ be vertices of $G$. Suppose that $x,y$ belong to the same
subset in the bipartitions of $G$ corresponding to
$C_1$ and $C_2$, respectively. Then $x,y$ belong to the same
subset in the bipartition of $G$ corresponding to $C_1\oplus C_2$.
\end{lem}

\proof
Immediate from Lemma~\ref{sum}.
\qed

\begin{lem}\label{linind}
Let $G$ be a connected graph and let $C_1,C_2,\ldots,C_m$ be  cocycles of $G$.
Suppose that for each $i$, there exists a pair of vertices  $x_i,y_i$ such
that $x_i,y_i$ belong to different subsets in the bipartition of $G$
corresponding to $C_i$
and to the same subset in the bipartition of $G$ corresponding to $C_j$
for all $j \neq i$ ($1\leq  j \leq m$).
Then  $C_1,C_2,\ldots,C_m$ are linearly independent.
\end{lem}

\proof
If not, then we may suppose without loss of generality that
$C_1=C_2\oplus C_3\oplus\ldots\oplus C_m$. This contradicts the fact that
$x_1,y_1$ belong to the different subsets
in the bipartition corresponding to $C_1$
and to the same subset in the bipartition corresponding to
$C_2\oplus C_3\oplus\ldots\oplus C_m$,
by Lemma~\ref{sameside}.
\qed

Let $G$ be a connected graph
and let $T$ be a tree on the same vertex set $V$ as $G$.
(We emphasize that $T$ need not be a subgraph of $G$.)
Each edge $e\in E(T)$ induces a bipartition of $V$
into nonempty subsets $X,Y$ given by the two components of $T-e$;
we define the {\em elementary cocycle of $G$ corresponding to $e$ and $T$}\/
to be the cocycle $E_G(X,Y)$.

\begin{lem}\label{basis}
Let $G$ be a connected graph with $n$ vertices,
and let $T$ be a tree on the same set of vertices
(not necessarily a subgraph of $G$).
For each edge $e_i\in E(T)$, let $C_i$ be the elementary
cocycle of $G$ corresponding to $e_i$ and $T$.
Then $\{C_1,C_2,\ldots,C_{n-1}\}$ is a basis for the  cocycle space of $G$.
\end{lem}

\proof
Using  Lemma~\ref{linind} (taking $x_i,y_i$ to be the end-vertices
of $e_i$) we deduce that $C_1,C_2,\ldots,C_{n-1}$ are linearly independent.
Since the dimension of the cocycle space of $G$ is $n-1$, they form a basis.
\qed

\begin{lem}\label{diffsidebasis}
Let $G$ be a connected graph with $n$ vertices,
let $\{C_1,C_2,\ldots,C_{n-1}\}$ be a basis for the cocycle space of $G$,
and let $x,y\in V(G)$ with $x \neq y$.
Then $x,y$ belong to different subsets in the bipartition corresponding
to $C_i$, for some $1\leq i\leq n-1$.
\end{lem}

\proof
Suppose not. Let $C$ be a cocycle in $G$ that separates $x$ and $y$
[for example, $E(\{x\}, \{x\}^c)$].
Since $\{C_1,C_2,\ldots,C_{n-1}\}$ is a basis for the  cocycle space of $G$,
$C$ is a linear combination of $C_1,C_2,\ldots,C_{n-1}$. This contradicts
Lemma~\ref{sameside}.
\qed

Now let $G$ be equipped with a family of
nonnegative real edge weights ${\bf w} = \{ w_e \} _{e \in E(G)}$.
As in \reff{def_weighted_flow}/\reff{def_weighted_maxmaxflow},
we let $\lambda_G(x,y;{\bf w})$ be the max flow from $x$ to $y$
with edge capacities ${\bf w}$,
and $\Lambda(G,{\bf w})$ the corresponding maxmaxflow.
As in \reff{def_weighted_cocycle},
we let $\widetilde{\Lambda}(G,{\bf w})$
be the minmax weight of the cocycles in a basis.
In order to prove the fundamental result (\ref{dlambdatwiddle}),
we shall need the following classic result on flows
(see \cite[Section 2.3]{Lovasz_86} for an excellent exposition):

\begin{theorem}[Gomory and Hu \protect\cite{GH}]  \label{gh}
Let $G$ be a connected graph equipped with nonnegative real edge weights
${\bf w} = \{w_e\} _{e \in E(G)}$.
Then there exists a tree $T$ with vertex set $V(T)=V(G) \equiv V$
(note that $T$ is not necessarily a subgraph of $G$!)
and a set ${\bf w}^T = \{w^T_e\} _{e \in E(T)}$
of nonnegative real edge weights such that
\begin{enumerate}
\item $\lambda_G(x,y;{\bf w})=\lambda_T(x,y;{\bf w}^T)$ for all $x,y\in V$
   ($x \neq y$), and
\item for each  $e=xy\in E(T)$, the elementary cocycle $C$ of $G$
corresponding to $e$ and $T$ is a minimum-weight edge cut
separating $x$ from $y$ in $G$, i.e.\
$\lambda_G(x,y;{\bf w})=\sum\limits_{f\in C}{w}_f$.
\end{enumerate}
\end{theorem}
We shall call any tree $T$ with the above properties a {\em Gomory--Hu tree}\/
for $(G, {\bf w})$;  it is in general nonunique.
Note that, for any given tree $T$, there is at most one choice of ${\bf w}^T$
that satisfies (a), namely for each edge $e = xy \in E(T)$
we must set $w^T_e = \lambda_G(x,y;{\bf w})$.
It can also be shown that if $T$ satisfies (b),
then this definition of ${\bf w}^T$ necessarily satisfies (a);
but we shall not need this fact.

If $T$ is a Gomory--Hu tree for $(G, {\bf w})$,
we define $\widehat{\Lambda}(G,{\bf w};T) = \max_{e \in E(T)} w^T_e$.
We claim that this value is independent of the choice of $T$,
and in fact we have:

\begin{thm}
\label{thdlambdatwiddle}
Let $G$ be a connected graph equipped with
nonnegative real edge weights ${\bf w} = \{ w_e \} _{e \in E(G)}$,
and let $T$ be a Gomory--Hu tree for $(G, {\bf w})$.
Then
\be
 \Lambda(G,{\bf w})  \;=\;
\widehat{\Lambda}(G,{\bf w};T)  \;=\;
\widetilde{\Lambda}(G,{\bf w})
    \;\le\; \Delta_2(G,{\bf w})  \;\le\; \Delta(G,{\bf w})
   \;.
\ee
In particular, the value of $\widehat{\Lambda}(G,{\bf w};T)$
is independent of the choice of $T$.
\end{thm}

\proof
The equality $\Lambda(G,{\bf w}) = \Lambda(T, {\bf w}^T)$
follows from Theorem~\ref{gh}(a),
and it is trivial to see that
$\Lambda(T, {\bf w}^T) = \max_{e \in E(T)} w^T_e$.
This proves that $\Lambda(G,{\bf w}) = \widehat{\Lambda}(G,{\bf w};T)$
and in particular that the latter quantity is independent of the choice of $T$.

The inequality $\Lambda(G,{\bf w}) \le \widetilde{\Lambda}(G,{\bf w})$
follows from Lemma~\ref{diffsidebasis}.

The inequality
$\widetilde{\Lambda}(G,{\bf w}) \le \widehat{\Lambda}(G,{\bf w};T)$
follows from Lemma~\ref{basis}
and Theorem~\ref{gh}(a,b).

There are easy elementary proofs of both
$\Lambda(G,{\bf w}) \le \Delta_2(G,{\bf w})$ and
$\widetilde{\Lambda}(G,{\bf w}) \le \Delta_2(G,{\bf w})$,
as noted in the Introduction.
\qed


\begin{corollary}
   \label{cordlambdatwiddle}
Let $G$ be a (not necessarily connected) graph equipped with
nonnegative real edge weights ${\bf w} = \{ w_e \} _{e \in E(G)}$.
Then
\be
 \Lambda(G,{\bf w})  \;=\;
\widetilde{\Lambda}(G,{\bf w})
    \;\le\; \Delta_2(G,{\bf w})  \;\le\; \Delta(G,{\bf w})
   \;.
\ee
\end{corollary}

\proof
If $G$ is disconnected, it suffices to apply Theorem~\ref{thdlambdatwiddle}
to each component of $G$.
\qed

Finally, we need to prove our claims that $\Lambda(G,{\bf w}) \ge D(G,{\bf w})$
and $\Lambda(G,{\bf w}) \ge$ \linebreak $\Delta_{n-1}(G,{\bf w})$.
We shall actually prove a slightly stronger result.
Define the {\em $k^{th}$ weighted degeneracy number}\/
\be
   D_k(G,{\bf w})   \;=\;   \max_{H \subseteq G}  \delta_k(H,{\bf w})
   \;,
\ee
where the max runs over all subgraphs $H$ of $G$,
and $\delta_k(H,{\bf w})$ denotes the $k$th smallest weighted degree of $H$:
\be
   \delta_k(H,{\bf w})   \;=\;
       \max\limits_{x_1,\ldots,x_{k-1} \in V(H)}
       \quad
       \min\limits_{x \in V(H) \backslash \{x_1,\ldots,x_{k-1}\}}
    d_H(x,{\bf w})
   \;.
\ee
Trivially we have
$D(G,{\bf w}) \equiv D_1(G,{\bf w}) \le D_2(G,{\bf w}) \le \ldots$
and $\delta_k(G,{\bf w}) \le D_k(G,{\bf w})$.
In particular,
\be
   D_2(G,{\bf w})   \;\ge\;   \max[D(G,{\bf w}), \, \delta_2(G,{\bf w})]
   \;.
\ee

\begin{proposition}
   \label{prop2.10}
Let $G$ be a graph with $n$ vertices ($n \ge 2$) equipped with
nonnegative real edge weights ${\bf w} = \{ w_e \} _{e \in E(G)}$.
Then
\be
   \Lambda(G,{\bf w})  \;\ge\;  D_2(G,{\bf w})
      \;\ge\;   \max[D(G,{\bf w}), \, \Delta_{n-1}(G,{\bf w})]
   \;.
\ee
\end{proposition}

\proof
Suppose first that $G$ is connected,
and let $T$ be a Gomory--Hu tree for $(G, {\bf w})$.
For any vertex $x$ of degree 1 in $T$,
let $e=xy$ be the unique incident edge in $T$;
then the elementary cocycle of $G$ corresponding to $e$ and $T$
is $E_G(\{x\}, \{x\}^c)$. Using Theorem~\ref{gh}(b), we have
$\Lambda(G,{\bf w})\geq \lambda_G(x,y;{\bf w})=d_G(x,{\bf w})$.
Since there are at least two such vertices $x$,
we have $\Lambda(G,{\bf w}) \ge \delta_2(G,{\bf w})$.

If $G$ is disconnected, we can apply the result just proven
to each component of $G$;
we conclude again that $\Lambda(G,{\bf w}) \ge \delta_2(G,{\bf w})$.

Now apply this result to each subgraph $H$ of $G$:
we conclude that $\Lambda(H,{\bf w}|_H) \ge \delta_2(H,{\bf w})$.
But $\Lambda(G,{\bf w}) \ge \Lambda(H,{\bf w}|_H)$
for every subgraph $H$ of $G$, so $\Lambda(G,{\bf w}) \ge D_2(G,{\bf w})$.
\qed

Let us now prove a few further general properties of maxmaxflow. Let
$G$ be a graph and $x\in V(G)$.
We say that $x$ is a {\em cut vertex}\/ of $G$
if $G\setminus x$ has more components than $G$.
We say that $G$ is {\em non-separable}\/
if $G$ is connected and has no cut vertices.\footnote{
   This concept is closely related to the more common notion of
   2-connectedness.
   A graph $G$ is {\em 2-connected}\/
   if $G$ has at least three vertices and $G \setminus x$ is connected
   for all $x \in V(G)$.
   Thus, a graph with at least three vertices is non-separable
   if and only if it is 2-connected.
   However, the graphs $K_1$ (a single vertex with no edges)
   and $K_2^{(m)}$ (a pair of vertices connected by $m$ parallel edges,
   with $m \ge 1$) are non-separable but not 2-connected.
}
A {\em block}\/ of $G$ is a maximal non-separable
subgraph of $G$. We first observe that maxmaxflow has a naturalness
property that maximum degree and $k$th-largest degree lack, namely,
it ``trivializes over blocks'':

\begin{proposition}
   \label{prop2.triv_blocks}
Let $G_1,\ldots,G_b$ be the blocks of $G$.
Then $\Lambda(G,{\bf w}) = \max\limits_{1 \le i \le b} \Lambda(G_i,{\bf w})$.
\end{proposition}

\proof
If $x$ and $y$ lie in the same block $G_i$,
then $\lambda_G(x,y;{\bf w}) = \lambda_{G_i}(x,y;{\bf w})$.
If $x$ and $y$ lie in the same component of $G$ but in different blocks,
then there exist cut vertices $v_1,\ldots,v_k$ of $G$
and blocks $G_{i_0},G_{i_1},\ldots,G_{i_k}$ of $G$
such that $x \in V(G_{i_0})$, $y \in V(G_{i_k})$,
$V(G_{i_{j-1}}) \cap V(G_{i_j}) = \{v_j\}$,
and every path from $x$ to $y$ passes through $v_1,\ldots,v_k$
in that order;
and in this case we have
\be
   \lambda_G(x,y;{\bf w})  \;=\;
   \min\Bigl[ \lambda_{G_0}(x,v_1;{\bf w}), \,
              \lambda_{G_1}(v_1,v_2;{\bf w}), \,
              \ldots, \,
              \lambda_{G_{k-1}}(v_{k-1},v_k;{\bf w}), \,
              \lambda_{G_k}(v_k,y;{\bf w})
       \Bigr]
   \;.
\ee
Finally, if $x$ and $y$ lie in different components of $G$,
then $\lambda_G(x,y;{\bf w}) = 0$.
\qed

It follows immediately from the definition of maxmaxflow that
for any pair of distinct vertices $x_1,x_2$ of $G$,
there exists a partition $V(G) = X_1 \cup X_2$ such that
$x_1 \in X_1$, $x_2 \in X_2$ and
$\sum_{e \in E(X_1,X_2)} w_e = \lambda_G(x_1,x_2;{\bf w})
 \le \Lambda(G,{\bf w})$.
We will need the following extension of this observation in Section
\ref{sec.block}.

\begin{proposition}
   \label{prop2.cuts}
Let $X\subseteq V(G)$ with $|X|\geq 2$. Then there exist
$x_1,x_2\in X$ and
disjoint $ X_1,X_2\subseteq V(G)$
such that  $X\cap X_i=\{x_i\}$ and
$\sum_{e \in E(X_i,X_i^c)} w_e \le \Lambda(G,{\bf w})$
for all $1\leq i\leq 2$.
\end{proposition}

\proof
We can assume without loss of generality that $G$ is connected.
Let $T$ be a Gomory--Hu tree for $(G,{\bf w})$,
and let $T'$ be the union of all the paths in $T$ connecting
pairs of vertices of $X$.
Let $x_1,x_2$ be distinct end-vertices of $T'$
and $e_1,e_2$ be the edges of  $T'$ incident with $x_1,x_2$
respectively.
For $i=1,2$, let $X_i$ be the vertex set of the component
of $T\setminus e_i$ which contains $x_i$. Then $X_1\cap X_2=\emptyset$.
Furthermore, $E(X_i,X_i^c)$ is the
elementary cocycle of $G$ corresponding to the edge $e_i$ of $T$;
so by Theorem~\ref{gh}(b) we have
$\sum_{e \in E(X_i,X_i^c)} w_e \le \Lambda(G,{\bf w})$
for $i=1,2$.
\qed

\bigskip

We conclude this section with a few examples of maxmaxflow calculations:

\bexam
 \label{example2.1}
Let $G$ be any forest, and let ${\bf w}$
be any set of nonnegative edge weights.
Then $\Lambda(G,{\bf w}) = \max_{e \in E(G)} w_e$.
Indeed, this elementary fact was already used in the proof of
Theorem~\ref{thdlambdatwiddle}.
\eexam

\bexam
 \label{example2.2}
Let $G$ be the cycle $C_n$, and let ${\bf w}$
be any set of nonnegative edge weights.  Then
\be
   \Lambda(G,{\bf w}) \;=\;
   \max\limits_{e \in E(G)} w_e \,+\, \min\limits_{e \in E(G)} w_e
   \;.
\ee
{\sc Proof.}
For any pair of distinct vertices $x,y \in V(G)$,
there are precisely two paths from $x$ to $y$,
and together they use all the edges of $G$.
Hence on one path the max flow is exactly $\min w_e$,
and on the other it is at most $\max w_e$.
So $\Lambda(G,{\bf w}) \le \min w_e + \max w_e$.
On the other hand, if we take $x,y$ to be the endpoints of the edge
with maximum weight, we obtain equality.
\eexam

\section{Some Preliminaries}   \label{sec.prelim}

Before turning to the counting of subgraphs of various classes,
let us make some brief prefatory observations.

\subsection{Pointwise bounds vs.\ generating-function bounds}

Let $(G,{\bf w})$ be a weighted graph,
and suppose that $a_m$ is the total weight
of all $m$-edge subgraphs of $G$ of some specified class.
In the following sections we shall prove upper bounds on $a_m$
of two different types:
\begin{itemize}
   \item ``Pointwise bound'':  $a_m \le C_m$
      for some specified constants $C_m$.
   \item ``Generating-function bound'':
      $\sum_{m=0}^\infty a_m z^m \le C(z)$
      for some specified function $C(z)$
      [which is allowed to take the value $+\infty$], for all $z \ge 0$.
\end{itemize}
{}From each type of bound we can deduce one of the other type:
$a_m \le C_m$ trivially implies
$\sum a_m z^m \le \sum C_m z^m$;
and $\sum a_m z^m \le C(z)$ implies
$a_m \le \inf_{z \ge 0} z^{-m} C(z)$.
But we lose something in each such passage:  for example,
\be
   \sum\limits_{m=0}^\infty a_m \le 1  \;\;\Longrightarrow\;\;
   a_m \le 1 \;\;\Longrightarrow\;\;
   \sum\limits_{m=0}^\infty a_m z^m \le {1 \over 1-z}
        \;\;\Longrightarrow\;\;
   a_m \le {(m+1)^{m+1} \over m^m} \;\;\Longrightarrow\;\; \cdots
\ee
So it is important to choose the best type of bound in each problem.
Roughly speaking, the pointwise bound is better in those cases
where the bound can be saturated simultaneously for all
(or arbitrarily many) values of $m$.
By contrast, in cases where the pointwise bound can be saturated
only for one value of $m$ at a time,
it may be possible to prove a stronger bound of generating-function type.
These remarks will be illustrated repeatedly throughout the
following sections.

\subsection{Convex hulls in graphs}

Let $H$ be any graph, and let $X$ be a nonempty subset of $V(H)$.
We define the {\em convex hull of $X$ in $H$}\/,
denoted ${\rm conv}(X,H)$,
to be the union of all paths in $H$ connecting
any pair of vertices $x_1,x_2 \in X$
(including paths of length 0 from a vertex $x \in X$ to itself).
Thus, ${\rm conv}(X,H)$ is a subgraph of $H$
whose vertex set contains $X$.
The following properties are elementary consequences of this
definition:
\begin{quote}
\begin{itemize}
  \item[(Conv1)] If $H$ is connected, then ${\rm conv}(X,H)$ is connected.
  \item[(Conv2)] Any vertex of degree 0 or 1 in ${\rm conv}(X,H)$
      must belong to $X$.  Moreover, a vertex $x \in X$ is of degree 0
      (i.e.\ isolated) in ${\rm conv}(X,H)$ if and only if
      it is the only element of $X$ in its component of $H$.
  \item[(Conv3)] If $H$ is a tree, then ${\rm conv}(X,H)$
      is the smallest subtree of $H$ containing all the vertices of $X$.
  \item[(Conv4)] If $H$ is non-separable and $|X| \ge 2$,
      then ${\rm conv}(X,H) = H$.
      [This is because, for any pair of distinct vertices $x_1,x_2$ of $H$,
       every edge of $H$ lies on some path from $x_1$ to $x_2$.]
\end{itemize}
\end{quote}
Property (Conv4) can be generalized as follows.  Let $H$ be a graph
and $x\in V(H)$. We say that $x$ is an {\em internal vertex}\/ of $H$
if $x$ is not a cut vertex of $H$. An {\em end block}\/ of $H$ is a
block that contains exactly one cut vertex of $H$.

\begin{lemma}
   \label{lemma.convexhull}
Let $H$ be connected, and let $X \subseteq V(H)$ with $|X| \ge 2$.
Then:
\begin{itemize}
   \item[(a)]  ${\rm conv}(X,H)$ is a connected union of blocks of $H$.
   \item[(b)]  Each end block of ${\rm conv}(X,H)$
      has an internal vertex belonging to $X$.
\end{itemize}
\end{lemma}

\proof This follows easily from the definition of ${\rm conv}(X,H)$ and
property (Conv4). \qed

\noindent
We have already used convex hulls in the proof of
Proposition~\ref{prop2.cuts},
and they will play an important role
in our treatment of trees and block trees
(Sections~\ref{sec.trees} and \ref{sec.block}).

\section{Counting Walks and Paths}  \label{sec.walks}

Let $G$ be a graph equipped with
nonnegative real edge weights ${\bf w} = \{ w_e \} _{e \in E(G)}$.
In this section (as well as in the following ones)
we shall write $\Delta$, $\Delta_2$, $\Lambda$, \ldots\
as a shorthand for
$\Delta(G,{\bf w})$, $\Delta_2(G,{\bf w})$, $\Lambda(G,{\bf w})$, \ldots;
the underlying graph $G$ and its edge weights ${\bf w}$
will always be understood.
Similarly, we shall write $\lambda(x,y)$
as a shorthand for $\lambda_G(x,y;{\bf w})$.

For $x,y \in V(G)$ and $m \ge 0$,
let $\scrw_m(x,y)$ be the set of $m$-step walks from $x$ to $y$,
i.e.\ sequences $\omega = x_0 e_1 x_1 e_2 x_2 \cdots x_{m-1} e_m x_m$
with $x_0=x$ and $x_m=y$ such that each $e_i$ is an edge $x_{i-1} x_i$.
We then define the following subsets of $\scrw_m(x,y)$:
\begin{itemize}
   \item $\scrw^{\rm FPW}_m(x,y)$ consists of {\em first-passage walks}\/,
       i.e.\ walks such that $x_i \neq y$ for $i < m$;
   \item $\scrw^{\rm SAW}_m(x,y)$ consists of {\em self-avoiding walks}\/
       (otherwise known as {\em paths}\/),
       i.e.\ walks such that $x_i \neq x_j$ for $i \neq j$.
       [We remark that a walk $\omega$ is self-avoiding if and only if
        each initial segment of $\omega$ is a first-passage walk.]
\end{itemize}
Clearly $\scrw^{\rm SAW}_m(x,y) \subseteq \scrw^{\rm FPW}_m(x,y)
\subseteq \scrw_m(x,y)$.

More generally, let $x \in V(G)$ and $Y \subseteq V(G)$.
We let $\scrw_m(x,Y) = \bigcup_{y \in Y} \scrw_m(x,y)$
be the set of $m$-step walks from $x$ to $Y$.
We then define the following subsets of $\scrw_m(x,Y)$:
\begin{itemize}
   \item $\scrw^{\rm FPW}_m(x,Y)$ consists of
       {\em first-passage walks to the set $Y$}\/,
       i.e.\ walks such that $x_i \notin Y$ for $i < m$;
   \item $\scrw^{\rm FPSAW}_m(x,Y)$ consists of
       {\em first-passage self-avoiding walks to the set $Y$}\/,
       i.e.\ walks such that $x_i \notin Y$ for $i < m$
       and $x_i \neq x_j$ for $i \neq j$.
\end{itemize}
Then $\scrw^{\rm FPSAW}_m(x,Y) \subseteq \scrw^{\rm FPW}_m(x,Y)
\subseteq \scrw_m(x,Y)$,
and in fact
$\scrw^{\rm FPSAW}_m(x,Y) =$ \linebreak $\scrw^{\rm FPW}_m(x,Y) \cap
 \bigcup\limits_{y \in Y} \scrw^{\rm SAW}_m(x,y)$.

We define the {\em weight} of a walk
$\omega = x_0 e_1 x_1 e_2 x_2 \cdots x_{m-1} e_m x_m$
to be the product of its edge weights:
\be
   w(x_0 e_1 x_1 e_2 x_2 \cdots x_{m-1} e_m x_m)   \;=\;
      \prod\limits_{i=1}^m w_{e_i}
   \;.
\ee
Finally, we define the weighted counts
\begin{eqnarray}
   w_m(x,y)   & = &   \sum\limits_{\omega \in \scrw_m(x,y)}  w(\omega)
   \\[2mm]
   w_m(x,Y)   & = &   \sum\limits_{\omega \in \scrw_m(x,Y)}  w(\omega)
\end{eqnarray}
and likewise for $w^{\rm FPW}$, $w^{\rm SAW}$ and $w^{\rm FPSAW}$.

The weighted count of walks
with a fixed initial vertex and an arbitrary final vertex
can trivially be bounded in terms of maximum weighted degree:
\begin{proposition}
   \label{prop4.1}
For all $x \in V(G)$ and all $m \ge 0$, we have
\be
   \sum_{y \in V(G)} w_m(x,y)  \;\le\;   \Delta^m   \;.
 \label{eq.walkbound}
\ee
\end{proposition}
In fact, when $(G,{\bf w})$ is $\Delta$-regular
[i.e.\ $d_G(x,{\bf w}) = \Delta$ for all $x \in V(G)$],
the bound \reff{eq.walkbound} is sharp simultaneously
for all $x \in V(G)$ and all $m \ge 0$.


If we restrict to first-passage walks to a fixed target set $Y$,
we can obtain a ``generating-function'' bound that is
sharper than \reff{eq.walkbound}:

\begin{prop}
   \label{prop4.2}
For all $x \in V(G)$ and all $Y \subseteq V(G)$, we have
\be
   \sum\limits_{m=0}^\infty \Delta^{-m} w_m^{\rm FPW}(x,Y)  \;\le\;  1
   \;.
 \label{eq.prop4.2}
\ee
\end{prop}

\proof
Define a sub-Markov chain\footnote{
   A {\em sub-Markov chain}\/ on a finite or countably infinite
   state space $X$ is defined by a transition kernel
   $\{ p(x \to x') \}_{x,x' \in X}$ satisfying
   $p(x \to x') \ge 0$ and $\sum_{x' \in X} p(x \to x') \le 1$.
   It induces a Markov chain on the state space $X \cup \{\infty\}$
   (where the new state $\infty \notin X$ is called the ``cemetery'')
   by defining $p(x \to\infty) = 1 - \sum_{x' \in X} p(x \to x')$,
   $p(\infty \to x) = 0$ for all $x \in X$,
   and $p(\infty \to \infty) = 1$.
}
with transition probabilities
$p(x \to x') = \Delta^{-1} \! \sum\limits_{e = xx'} w_e$
(here the sum runs over
all edges $e=xx'$ in $G$). Then $\sum_{m=0}^\infty \Delta^{-m}
w_m^{\rm FPW}(x,Y)$ is the probability that $Y$ is eventually hit,
starting at $x$. \qed

\medskip\noindent
{\bf Remarks.}
1.  If $G$ is connected (with weights $w_e > 0$)
and $(G,{\bf w})$ is $\Delta$-regular,
then \reff{eq.prop4.2} is {\em equality}\/
for all $x \in V(G)$ and all nonempty $Y \subseteq V(G)$.
This is because $p$ is a finite irreducible Markov chain,
so that $Y$ is hit with probability 1.
Indeed, it suffices to have $d_G(x,{\bf w}) = \Delta$
for all $x \in V(G) \setminus Y$;
the degree at vertices of $Y$ is irrelevant.

2.  If $x \in Y$ then $w_0^{\rm FPW}(x,Y)=1$ and $w_m^{\rm
FPW}(x,Y)=0$ for $m\geq 1$ so equality holds in \reff{eq.prop4.2}.
On the other hand, if $x \notin Y$, then we can improve the upper
bound in \reff{eq.prop4.2} from 1 to $d_G(x,{\bf w})/\Delta$, since
this is the probability for the Markov chain to survive the first
step.

3.  As just mentioned, \reff{eq.prop4.2} is the best-possible upper bound
for first-passage walks to a set $Y$.
However, one might ask whether a sharper bound is possible
by restricting attention to first-passage {\em self-avoiding}\/ walks.
The answer is no, at least if one considers a general set $Y$,
as the following examples show:

\bexam
  \label{example4.1}
Let $G$ be the star $K_{1,r}$,
let $x$ be the central vertex,
and let $Y$ be the remaining vertices.
Set all edge weights $w_e = \Delta/r$.
Then $w_1^{\rm FPSAW}(x,Y) = \Delta$
and $w_m^{\rm FPSAW}(x,Y) = 0$ for $m \neq 1$,
so \reff{eq.prop4.2} is sharp.
\eexam

\bexam
  \label{example4.2}
More generally, let ${\bf T}_r$ be the infinite $r$-regular tree,
let $x$ be any vertex of ${\bf T}_r$, and fix $n \ge 1$.
Let $G$ be the subtree of ${\bf T}_r$
induced by the vertices lying at a distance at most $n$ from $x$,
and let $Y$ be the set of vertices lying at a distance exactly $n$ from $x$.
Set all edge weights $w_e = \Delta/r$.
Then
\be
   w_m^{\rm FPSAW}(x,Y) \;=\;
   \cases{ r(r-1)^{m-1} (\Delta/r)^m   & if $m = n$  \cr
           0                           & otherwise   \cr
         }
\ee
so that
\be
   \sum_{m=0}^\infty \Delta^{-m} w_m^{\rm FPSAW}(x,Y) \;=\;
   \left({r-1 \over r}\right) ^{\! n-1}
   \;,
\ee
which is sharp in the limit $r \to \infty$ at any fixed $n$.
(It is also sharp for any $r$ when $n=1$, which is the case $G=K_{1,r}$.)
\eexam

\bigskip

For comparison with Proposition~\ref{prop4.3} below,
it is perhaps illuminating to rephrase the proof of
Proposition~\ref{prop4.2} as an induction:

\par\bigskip\noindent{\sc Second proof of Proposition~\ref{prop4.2}.\ }
We will prove inductively that
$$\sum_{m=0}^M \Delta^{-m} w_m^{\rm FPW}(x,Y) \le 1$$
for all $M \ge 0$ and all $x \in V(G)$.
This clearly holds for $M=0$.
Moreover, it holds for all $M$ when
$x \in Y$, since $w_m^{\rm FPW}(x,Y) = \delta_{m0}$.
So suppose that $x \notin Y$.
Then $w_0^{\rm FPW}(x,Y) = 0$ and for $m \ge 1$ we have
\be
   w_m^{\rm FPW}(x,Y)   \;=\;
   \sum\limits_{x' \in V(G)} \left( \sum\limits_{e = xx'} w_e \right) \,
      w_{m-1}^{\rm FPW}(x',Y)
   \;.
\ee
Thus,
for $M \ge 1$, we have
\begin{eqnarray}
   \sum\limits_{m=0}^M \Delta^{-m} w_m^{\rm FPW}(x,Y)
      & = &
   \sum\limits_{m=1}^M \Delta^{-m}
      \sum\limits_{x' \in V(G)} \left( \sum\limits_{e= xx'} w_e \right) \,
         w_{m-1}^{\rm FPW}(x',Y)
   \nonumber \\[2mm]
      & = &
   \Delta^{-1} \sum\limits_{x' \in V(G)}
                    \left( \sum\limits_{e = xx'} w_e \right) \,
               \sum\limits_{j=0}^{M-1} \Delta^{-j} w_{j}^{\rm FPW}(x',Y)
   \nonumber \\[2mm]
      & \le &   1  \quad\hbox{by the inductive hypothesis.}
\end{eqnarray}
\qed


There is no hope of obtaining a $|Y|$-independent bound on
$w_m^{\rm FPW}(x,Y)$ --- or even on $w_m^{\rm FPSAW}(x,Y)$ ---
in terms of second-largest degree (much less in terms of maxmaxflow).
Indeed, let $G$ be the star $K_{1,r}$,
so that $\Delta = r$ and $\Delta_2 = \Lambda = 1$;
let $x$ be the central vertex and $Y$ the remaining vertices;
then $w_1^{\rm FPSAW}(x,Y) = r$, which is unbounded.
Nor can we bound $w_m^{\rm FPW}(x,y)$:
considering again the star and taking $y$ to be any vertex other than $x$,
we have $w_3^{\rm FPW}(x,y) = r-1$, which is again unbounded.
Nevertheless, an analogue of Proposition~\ref{prop4.2}
does hold for maxmaxflow if we restrict ourselves to self-avoiding walks
and to the case $Y = \{y\}$:

\begin{prop}
   \label{prop4.3}
Define
\be
   F(x,y)  \;=\;  \cases{ \lambda(x,y)/\Lambda   & if $x \neq y$  \cr
                          1                      & if $x=y$  \cr
                        }
\ee
Then, for all $x,y \in V(G)$, we have
\be
   \sum\limits_{m=0}^\infty \Lambda^{-m} w_m^{\rm SAW}(x,y)
   \;\le\;  F(x,y)  \;\le\; 1   \;.
 \label{eq.prop4.3}
\ee
\end{prop}

\proof
We will prove inductively that
$\sum_{m=0}^M \Lambda^{-m} w_m^{\rm SAW}(x,y) \le F(x,y)$
for all $M \ge 0$ and all $x,y \in V(G)$.
This clearly holds for $M=0$, since $w_0^{\rm SAW}(x,y) = \delta_{xy}$.
Moreover, it holds for all $M$ when $x=y$,
since $w_m^{\rm SAW}(x,x) = \delta_{m0}$.
So let $M \ge 1$ and $x \neq y$.
Let $C=E(X,Y)$ be a cocycle in $G$ with $x \in X$, $y \in Y$
and $\sum_{e \in C} w_e = \lambda(x,y)$.
Let $e = uv$ be an edge in $C$ with $u \in X$ and $v \in Y$,
and let $w_m^{\rm SAW}(x,y,e)$ be the weighted sum over $m$-step paths
from $x$ to $y$ that use $e$ as their first edge from $X$ to $Y$.
Then
\be
   w_m^{\rm SAW}(x,y,e)  \;\le\;
      w_e \sum\limits_{i=0}^{m-1} w_i^{\rm SAW}(x,u) w_{m-1-i}^{\rm SAW}(v,y)
   \;.
\ee
So
\begin{eqnarray}
   \sum\limits_{m=0}^M \Lambda^{-m} w_m^{\rm SAW}(x,y,e)
      & \le &
   \Lambda^{-1} w_e \sum\limits_{m=0}^M \sum\limits_{i=0}^{m-1}
     \Lambda^{-i} w_i^{\rm SAW}(x,u) \Lambda^{-(m-1-i)} w_{m-1-i}^{\rm SAW}(v,y)
   \nonumber \\[2mm]
      & = &
   \Lambda^{-1} w_e \sum\limits_{i=0}^{M-1} \Lambda^{-i} w_i^{\rm SAW}(x,u)
            \sum\limits_{j=0}^{M-1-i} \Lambda^{-j} w_j^{\rm SAW}(v,y)
   \nonumber \\[2mm]
      & \le &
   \Lambda^{-1} w_e \quad\hbox{by the inductive hypothesis}
\end{eqnarray}
(note that we may have here used the $x=y$ case of the inductive hypothesis,
 in case $x=u$ or $y=v$ or both).
Therefore
\begin{eqnarray}
   \sum\limits_{m=0}^M \Lambda^{-m} w_m^{\rm SAW}(x,y)
      & = &
   \sum\limits_{e \in C} \sum\limits_{m=0}^M \Lambda^{-m} w_m^{\rm SAW}(x,y,e)
   \nonumber \\[2mm]
      & \le &
   \sum\limits_{e \in C} \Lambda^{-1} w_e
   \nonumber \\[2mm]
      & \le &   \lambda(x,y)/\Lambda   \;.
 \label{eq.prop4.3.proof}
\end{eqnarray}
\qed

\begin{cor}
   \label{cor4.4}
$w_m^{\rm SAW}(x,y) \le \Lambda^m F(x,y)$
for all $m \ge 0$ and all $x,y \in V(G)$.
\end{cor}

\begin{cor}
   \label{cor4.5}
For $0 \le \zeta \le \Lambda^{-1}$, we have
$\sum\limits_{m=0}^\infty \zeta^m w_m^{\rm SAW}(x,y)
 \le (\zeta\Lambda)^{{\rm dist}(x,y)} F(x,y)$
where ${\rm dist}(x,y)$ is the length of the shortest path in $G$
from $x$ to $y$.
(Of course, if there is no such path,
 then $w_m^{\rm SAW}(x,y) = 0$ for all $m$.)
\end{cor}

\bexam
  \label{example4.3}
Let $G$ be any tree, let $w_e = w$ for all edges $e$,
and let $x,y$ be any two vertices of $G$.
Then $\lambda(x,y) = \Lambda = w$ and
\be
   w_m^{\rm SAW}(x,y) \;=\;
   \cases{ w^m   & if $m = {\rm dist}(x,y)$  \cr
           0     & otherwise \cr
         }
\ee
so that \reff{eq.prop4.3} is equality for all $x,y$.
Note that the generating-function bound \reff{eq.prop4.3}
is here much sharper than the pointwise bound of Corollary~\ref{cor4.4}:
this reflects the fact that
the latter can be saturated here for only one value of $m$
(for any given pair $x,y$).
This is the prototypical situation in which we shall seek
generating-function bounds.
Indeed, the inductive ``cutting'' argument used in the proof
of Proposition~\ref{prop4.3} seems to work {\em only}\/
for the generating-function bound;
we do not know of any way of proving Corollary~\ref{cor4.4}
without proving the stronger bound \reff{eq.prop4.3}.
\eexam

\bexam
  \label{example4.4}
Let $G$ be the complete graph $K_n$,
and let all edge weights $w_e$ equal $\Delta/(n-1)$,
so that both the maximum weighted degree and the maxmaxflow equal $\Delta$.
If $x,y$ are any two distinct vertices, we have
\be
   w_m^{\rm SAW}(x,y)  \;=\;
   {(n-2)(n-3) \cdots (n-m)  \over  (n-1)^m} \,  \Delta^m
\ee
if $1 \le m \le n-1$, and 0 otherwise.
Approximating Riemann sums by integrals, we find
\be
   \sum_{m=0}^\infty \Delta^{-m} w_m^{\rm SAW}(x,y)
   \;\approx\;
   \sqrt{{\pi \over 2n}}
\ee
as $n\to\infty$, so Proposition~\ref{prop4.3} is far from sharp
in this limit.  Nevertheless, we see that the exponential growth rate
$w_m^{\rm SAW}(x,y) \sim \Delta^m$
in Proposition~\ref{prop4.3} and Corollary~\ref{cor4.4}
cannot be improved.
Indeed, if for each fixed $m$ we choose $n$ so as to maximize
$w_m^{\rm SAW}(x,y)/\Delta^m$,
we find that the maximum is achieved at $n \approx m^2/2$
and that the maximum value is $\approx e/(2m^2)$.
\eexam

\bexam
  \label{example4.5}
Let $G$ be the generalized theta graph $\Theta_{1,2,\ldots,r}$ ($r
\ge 2$), consisting of a pair of endvertices $a,b$ joined by $r$
internally disjoint paths of lengths $1,2,\ldots,r$. On each path
let one edge have weight $w$ and the other edges have weight 1. Then
the maxmaxflow is \be
   \Lambda  \;=\;
   \cases{1+w      & if $0 \le w \le 1/(r-1)$   \cr
          rw       & if $1/(r-1) \le w \le 1$   \cr
          r-1+w    & if $w \ge 1$               \cr
         }
\ee while \be
   w_m^{\rm SAW}(a,b)  \;=\;
   \cases{w   & if $1 \le m \le r$  \cr
          0   & otherwise           \cr
         }
\ee
In particular, at $w=1/(r-1)$ we have

\be
   \sum_{m=0}^\infty \Lambda^{-m} w_m^{\rm SAW}(a,b)  \;=\;
   1 \,-\, \left(1 - {1 \over r}\right)^{\! r}
   \;,
\ee
which decreases to $1-1/e \approx 0.632121$ as $r \to\infty$.
So, although the bound of Proposition~\ref{prop4.3}
is not sharp in this case,
it does at least come within a constant factor of being sharp
in a situation where the maximum contribution from a single value of $m$
goes to zero
(the opposite extreme from Examples~\thesection.\ref{example4.1} and
 \thesection.\ref{example4.2}).
\eexam

\section{Counting Trees and Forests}   \label{sec.trees}

Let us now extend Propositions~\ref{prop4.2} and \ref{prop4.3}
from paths to trees and forests.
In Section~\ref{sec.trees.1} we consider classes $\scrt_m(X)$ of trees
and $\scrf_m(X,Y)$ of forests.
In Section~\ref{sec.trees.2} we consider
a larger class $\scrh_m(X)$ of forests.

\subsection[The classes $\scrt_m(X)$ and $\scrf_m(X,Y)$]{
   The classes \mbox{\boldmath$\scrt_m(X)$} and \mbox{\boldmath$\scrf_m(X,Y)$}}

  \label{sec.trees.1}

For $F$ a forest in $G$, let $L(F)$ denote the set of
vertices of degree 0 or 1 in $F$
(also called {\em leaves}\/ or {\em end-vertices}\/ of $F$).

For any nonempty $X \subseteq V(G)$, let $\scrt_m(X)$ be
the set of all $m$-edge trees $T$ in $G$ such that
$L(T)\subseteq X \subseteq V(T)$.
Heuristically, $\scrt_m(X)$ consists of trees whose leaves
are ``tied down'' on the set $X$.
Note the following special cases:
\begin{itemize}
   \item If $X = \{x\}$, then $\scrt_0(\{x\})$ has as its single element
the edgeless graph with vertex set $\{x\}$,
and $\scrt_m(\{x\}) = \emptyset$ for $m \ge 1$.
   \item $\scrt_m(\{x,y\}) \simeq \scrw_m^{\rm SAW}(x,y)$
under the natural identification of paths with their induced subgraphs;
this holds both for $x \neq y$ and for $x=y$.
   \item $\scrt_0(X) = \emptyset$ for $|X| \ge 2$.
\end{itemize}

For any $X,Y\subseteq V(G)$ with $Y \neq \emptyset$, let
$\scrf_m(X,Y)$ be the set of all $m$-edge forests $F$ in $G$ such
that
\begin{quote}
\begin{itemize}
   \item[(F1)] $L(F)\subseteq X\cup Y\subseteq V(F)$, and
   \item[(F2)] each component of $F$ contains exactly one vertex of $Y$.
\end{itemize}
\end{quote}
[Note that $\scrf_m(X,Y) = \scrf_m(X \setminus Y,Y)$,
so we can assume without loss of generality, if desired, that $X
\cap Y = \emptyset$.]
Heuristically, $\scrf_m(X,Y)$ consists of forests whose leaves
are ``tied down'' on the set $X\cup Y$ and whose components are ``tied down''
on single elements of the set $Y$.
We have the following special cases:
\begin{itemize}
   \item  If $X = \emptyset$ (or more generally if $X \subseteq Y$),
then $\scrf_0(X,Y)$ has as its single element the edgeless graph with
vertex set $Y$, and $\scrf_m(X,Y) = \emptyset$ for $m \ge 1$.
   \item  If $X = \{x\}$, then each $F \in \scrf_m(\{x\},Y)$
is the disjoint union of a path $P \in \scrw_m^{\rm FPSAW}(x,Y)$
[under the natural identification of paths with their induced
subgraphs] and the collection $Y \setminus V(P)$ of isolated
vertices; this holds both for $x \notin Y$ and for $x\in Y$. We
express this isomorphism loosely by writing $\scrf_m(\{x\},Y) \simeq
\scrw_m^{\rm FPSAW}(x,Y)$.
   \item  For $Y = \{y\}$, we have $\scrf_m(X,\{y\}) = \scrt_m(X \cup \{y\})$.
\end{itemize}

For $H$ a subgraph of $G$, set $w(H)=\prod_{e\in E(H)} w_e$.
[Note that if $E(H)=\emptyset$, then $w(H)=1$.]
Define the weighted counts
\begin{eqnarray}
   t_m(X)    & = &   \sum\limits_{T \in \scrt_m(X)}  w(T)    \\[2mm]
   f_m(X,Y)  & = &   \sum\limits_{F \in \scrf_m(X,Y)}  w(F)
\end{eqnarray}
We will obtain two ``generating-function'' bounds on $f_m(X,Y)$: a
sharp bound in terms of the maximum weighted degree $\Delta$, and a
slightly weaker bound in terms of the maxmaxflow $\Lambda$.

\begin{proposition}
  \label{prop.forests.Delta}
For all $X,Y\subseteq V(G)$ with $Y \neq \emptyset$, we have
\be
   \sum\limits_{m=0}^\infty \Delta^{-m} f_m(X,Y)
   \;\le\; 1   \;.
\ee
\end{proposition}

\begin{proposition}
  \label{prop.forests.Lambda}
For all $X,Y\subseteq V(G)$ with $Y \neq \emptyset$, we have
\be
   \sum\limits_{m=0}^\infty (m+|Y|)^{-(|X|-1)} \Lambda^{-m} f_m(X,Y)
   \;\le\; |Y|   \;.
\ee
\end{proposition}

Using the identity $\scrf_m(X,\{y\}) = \scrt_m(X \cup \{y\})$,
we immediately obtain the following corollaries for trees:

\begin{corollary}
  \label{cor.trees.Delta}
For all nonempty $X \subseteq V(G)$, we have
\be
   \sum\limits_{m=0}^\infty \Delta^{-m} t_m(X) \;\le\; 1   \;.
\ee
\end{corollary}

\begin{corollary}
  \label{cor.trees.Lambda}
For all nonempty $X \subseteq V(G)$, we have
\be
   \sum\limits_{m=0}^\infty (m+1)^{-(|X|-2)} \Lambda^{-m} t_m(X)
   \;\le\; 1   \;.
\ee
\end{corollary}

One can envisage two different approaches to proving these bounds
for forests and trees. One is to {\em imitate}\/ the ``cutting''
argument employed for walks in the proofs of
Propositions~\ref{prop4.2} and \ref{prop4.3}. The other is to {\em
exploit the result}\/ obtained for walks in
Propositions~\ref{prop4.2} and \ref{prop4.3}, by reducing forests
inductively to paths. The latter technique seems to be more appropriate.

In the proof of Proposition~\ref{prop.forests.Delta}
we shall make use of the ``point-to-set'' bound of Proposition~\ref{prop4.2}.
We shall also use the following elementary lemma
which  splits a forest with $k$ end-vertices into  a forest with
$k-1$ end-vertices
and a path:

\begin{lemma}
  \label{lemma.decompose}
Let $G$ be a graph, let $X,Y\subseteq V(G)$ with $Y \neq \emptyset$,
let $x\in X \setminus Y$, and let $F \in \scrf_m(X,Y)$.
Let $F_1$ be the convex hull of $(X \setminus x)\cup Y$ in $F$,
and let $P$ be the unique path in $F$ from $x$ to $V(F_1)$.
Then $F$ is the edge-disjoint union of $F_1$ and $P$;
and for some $i$ $(0 \le i \le m)$
we have $F_1 \in {\scrf}_i(X \setminus x,Y)$
and $P \in \scrw_{m-i}^{\rm FPSAW}(x,V(F_1))$.
Moreover, the map $F \mapsto (F_1,P)$ is a bijection
from $\scrf_m(X,Y)$ onto
$\bigcup\limits_{i=0}^m
 {\scrf}_i(X \setminus x,Y) \times \scrw_{m-i}^{\rm FPSAW}(x,V(F_1))$.
\end{lemma}

\noindent The proof is a straightforward exercise: let us simply
observe that because $x \notin Y$, the component of $F$ containing
$x$ must contain a  vertex in $Y$, which guarantees that the path
$P$ exists; and $P$ is unique because $F$ has no cycles.

\proofof{Proposition~\ref{prop.forests.Delta}}
As noted above, we can assume without loss of generality that
$X \cap Y = \emptyset$.
We use induction on $k=|X|$.
For $k=0$ the result is trivial.
For $k=1$ the result follows immediately from
Proposition~\ref{prop4.2},
since ${\scrf}_m(\{x\},Y) \simeq \scrw_m^{\rm FPSAW}(x,Y)
       \subseteq \scrw_m^{\rm FPW}(x,Y)$.
So suppose $k\geq 2$.
Let $x$ be any vertex in $X$;  by assumption $x \notin Y$.
By Lemma~\ref{lemma.decompose},
given any $F \in {\scrf}_m(X,Y)$,
we can decompose $F$ into a forest $F_1 \in {\scrf}_i(X \setminus x,Y)$
and a path $P \in \scrw_{m-i}^{\rm FPSAW}(x,V(F_1))$;
and each such pair $(F_1,P)$ arises from a unique $F$ (namely, $F_1 \cup P$).
Since $w(F)=w(F_1)w(P)$, we have
\begin{eqnarray}
& & \hspace*{-5mm}
 \sum_{m=0}^\infty \Delta^{-m} \!\!
 \sum_{F \in {\scrf}_m(X,Y)}   w(F)
  \nonumber \\
&&  = \sum_{m=0}^\infty
    \sum_{i=0}^m \Delta^{-i} \!\!
        \sum_{F_1\in {\scrf}_i(X \setminus x,Y)}
           w(F_1) \,
        \Delta^{-(m-i)} \!\!\!
    \sum_{P\in \scrw_{m-i}^{\rm FPSAW}(x,V(F_1))}
       w(P)
   \nonumber \\
&&  = \sum_{i=0}^\infty  \Delta^{-i} \!\!
      \sum_{F_1\in {\scrf}_i(X \setminus x,Y)}
         w(F_1)
      \sum_{j=0}^\infty  \Delta^{-j} \!\!\!
      \sum_{P\in \scrw_j^{\rm FPSAW}(x,V(F_1))}
     w(P)   \;.
\end{eqnarray}
Now, for each fixed $F_1$ we have
$\sum_{j=0}^\infty \Delta^{-j}
 \sum_{P \in \scrw_j^{\rm FPSAW}(x,V(F_1))} w(P) \leq 1$
by Proposition~\ref{prop4.2}, so that
\be
 \sum_{m=0}^\infty \Delta^{-m} \!\!
 \sum_{F\in {\scrf}_m(X,Y)}
     w(F)
   \;\le\;
 \sum_{i=0}^\infty \Delta^{-i} \!\!
 \sum_{F_1\in {\scrf}_i(X \setminus x,Y)}
     w(F_1)
   \;\le\;  1
\ee
by the inductive hypothesis.
\qed

Examples~\ref{sec.walks}.\ref{example4.1} and \ref{sec.walks}.\ref{example4.2}
show that Proposition~\ref{prop.forests.Delta} is best possible
(at least for a general set $Y$), even when $|X|=1$.

\proofof{Proposition~\ref{prop.forests.Lambda}}
As before, we assume that $X \cap Y = \emptyset$
and we use induction on $k=|X|$.
If $k=0$ the result is trivial.
Suppose next that $k=1$ and $X=\{x\}$.
Since $f_m(\{x\},\{y\})=w_m^{\rm SAW}(x,y)$, we have
$\sum\limits_{m=0}^\infty \Lambda^{-m} f_m(\{x\},\{y\})
   \;\le\; 1$ for each $y\in Y$ by Proposition~\ref{prop4.3}. Thus
$\sum\limits_{m=0}^\infty \Lambda^{-m} f_m(\{x\},Y)
   \;\le\; |Y|$ and the proposition holds for $k=1$.

Suppose now that $k\geq 2$.
Let $x$ be any vertex in $X$;  by assumption $x \notin Y$.
By Lemma~\ref{lemma.decompose},
given any $F \in {\scrf}_m(X,Y)$,
we can decompose $F$ into a forest
$F_1 \in {\scrf}_i(X \setminus x,Y)$
and a path $P \in \scrw_{m-i}^{\rm FPSAW}(x,V(F_1))$;
and each such pair $(F_1,P)$ arises from a unique $F$.
Since $w(F)=w(F_1)w(P)$, we have
\begin{eqnarray}
& & \hspace*{-5mm}
 \sum_{m=0}^\infty (m+|Y|)^{-(k-1)}\Lambda^{-m} \!\!
 \sum_{F \in {\scrf}_m(X,Y)}   w(F)
  \nonumber \\
&&  = \sum_{m=0}^\infty (m+|Y|)^{-(k-1)}
    \sum_{i=0}^m \Lambda^{-i} \!\!
        \sum_{F_1\in {\scrf}_i(X \setminus x,Y)}
           w(F_1) \,
        \Lambda^{-(m-i)} \!\!\!
    \sum_{P\in \scrw_{m-i}^{\rm FPSAW}(x,V(F_1))}
       w(P)
   \nonumber \\
&&  = \sum_{i=0}^\infty  \Lambda^{-i} \!\!
      \sum_{F_1\in {\scrf}_i(X \setminus x,Y)}
         w(F_1)
      \sum_{j=0}^\infty  \Lambda^{-j}  (i+j+|Y|)^{-(k-1)} \!\!\!
      \sum_{P\in \scrw_j^{\rm FPSAW}(x,V(F_1))}
         w(P)   \;.
   \nonumber \\
\end{eqnarray}
Now, for each fixed $F_1$ we have
$\sum_{j=0}^\infty \Lambda^{-j}
 \sum_{P \in \scrw_j^{\rm FPSAW}(x,V(F_1))} w(P) \leq |V(F_1)|$
by the base case $k=1$.
Since $|V(F_1)|= i+|Y| \le i+j+|Y|$, we have
\begin{eqnarray}
 \sum_{m=0}^\infty (m+|Y|)^{-(k-1)}
\Lambda^{-m} \!\!
 \sum_{F\in {\scrf}_m(X,Y)}
     w(F)
   & \le &
 \sum_{i=0}^\infty (i+|Y|)^{-(k-2)}
\Lambda^{-i} \!\!
 \sum_{F_1\in {\scrf}_i(X \setminus x,Y)}
     w(F_1)
    \nonumber \\
&\leq & 1
    \nonumber \\
\end{eqnarray}
by the inductive hypothesis.
\qed

When $|X|=1$, Proposition~\ref{prop.forests.Lambda} is in some sense
best possible. To see this, take $G$ to be a star, $X$ to be the
central vertex and $Y$ to be the end-vertices: this gives
$\Lambda=1$, $f_1(X,Y)=|Y|$ and $f_m(X,Y)=0$ for $m\neq 1$, so that
$\sum_{m=0}^\infty \Lambda^{-m} f_m(X,Y)=|Y|$. When $|X|=k\geq 2$,
by contrast, Proposition~\ref{prop.forests.Lambda} is perhaps not
best possible. If we take $G$ to be the disjoint union of $k$
isomorphic stars (again with central vertices in $X$ and
end-vertices in $Y$), we have $\sum_{m=0}^\infty
\Lambda^{-m}f_m(X,Y)=(|Y|/k)^k$. This shows that if there is a
universal upper bound on $\sum_{m=0}^\infty \Lambda^{-m}f_m(X,Y)$,
the right-hand side has to grow at least like $(|Y|/|X|)^{|X|}$. We
suspect the following conjectures are true:

\begin{conjecture}
   \label{conj.forests.Lambda}
For all $X,Y\subseteq V(G)$ with $Y \neq \emptyset$, we have
\be
   \sum_{m=0}^\infty \Lambda^{-m} f_m(X,Y)  \;\le\; |Y|^{|X|}  \;.
\ee
\end{conjecture}

\begin{conjecture}[a special case of Conjecture~\ref{conj.forests.Lambda}]
   \label{conj.trees.Lambda}
For all nonempty $X \subseteq V(G)$, we have
\be
   \sum_{m=0}^\infty \Lambda^{-m} t_m(X) \;\le\; 1   \;.
\ee
\end{conjecture}

\subsection[The class $\scrh_m(X)$]{The class \mbox{\boldmath$\scrh_m(X)$}}
   \label{sec.trees.2}

We conclude this section by discussing a larger class of forests.
For $X\subseteq V(G)$, let $\scrh_m(X)$ be
the set of all $m$-edge forests $F$ in $G$ such that
$L(F)\subseteq X\subseteq V(F)$.
For integers $p,r \ge 1$, let $\scrh_m(X,p)$ be
the set of all $m$-edge forests $F$ in $\scrh_m(X)$ such that
each component of $F$ contains at least $p$ vertices of $X$,
and let $\scrh_m(X,p,r)$ be
the set of all forests $F$ in $\scrh_m(X,p)$ such that
$F$ has precisely $r$ components. Put
\be
   h_m(X)  \;=\;   \sum\limits_{F \in \scrh_m(X)}  w(F)  \;,
\ee
and define $h_m(X,p)$ and $h_m(X,p,r)$ similarly.
Our next result uses Proposition~\ref{prop.forests.Delta} to bound
$h_m(X,p,r)$ in terms of $\Delta$:

\begin{proposition}
  \label{prop.h_m(X,p,r).Delta}
Let $X\subseteq V(G)$ where $|X|=k\geq rp$. Then
\be
   \sum\limits_{m=0}^\infty \Delta^{-m} h_m(X,p,r)
   \;\le\;    p^{-(r-1)}(k-rp+p)^{-1} {{k}\choose{r}}\;.
\ee
\end{proposition}
\proof
Choose $F\in \scrh_m(X,p,r)$.
Let
$\{X_1,X_2,\ldots,X_r\}$ be the partition of $X$ determined by the
components of $F$. Then $F\in \scrf_m(X,Y)$ for all sets $Y$ such that
$|Y\cap X_j|=1$ for all $1\leq j\leq r$.
Hence there are precisely $\prod_{j=1}^r|X_j|$ different
sets $Y\subseteq X$ such that $F\in \scrf_m(X,Y)$.
Since $|X_j|\geq p$ for all $j$ ($1\leq j\leq r$)
and $\sum_{j=1}^r |X_j| = k$,
it follows that $\prod_{j=1}^r|X_j|\geq p^{r-1}(k-rp+p)$.
Thus
\be
      \sum\limits_{F\in \scrh_m(X,p,r)} w(F)
      \;\le\;
p^{-(r-1)}(k-rp+p)^{-1}
   \sum\limits_{\begin{scarray}
             Y\subseteq X\\
                |Y|=r
        \end{scarray}}
   \sum\limits_{F\in \scrf_m(X,Y)}w(F).
\ee
Using Proposition~\ref{prop.forests.Delta} we deduce
\begin{eqnarray}
      \sum\limits_{m=0}^\infty
      \sum\limits_{F\in \scrh_m(X,p,r)} \Delta^{-m} w(F)
      &\le &
p^{-(r-1)}(k-rp+p)^{-1}
   \sum\limits_{m=0}^\infty  \Delta^{-m}
   \sum\limits_{\begin{scarray}
             Y\subseteq X\\
                |Y|=r
        \end{scarray}}
   \sum\limits_{F\in \scrf_m(X,Y)}w(F)    \nonumber \\
   &\le&
   p^{-(r-1)}(k-rp+p)^{-1}{{k}\choose{r}}\;.
\end{eqnarray}
\qed

Summing over $r$,  we obtain:
\begin{corollary}
  \label{prop.h_m(X,p).Delta}
Let $X\subseteq V(G)$ where $|X|=k\geq 1$.  Then
\begin{subeqnarray}
   \sum\limits_{m=0}^\infty \Delta^{-m} h_m(X,p)
   & \le &
   \sum\limits_{r=1}^{\lfloor k/p \rfloor} p^{-(r-1)}(k-rp+p)^{-1} {k\choose r}
      \slabel{eq.prop.h_m(X,p).Delta.a}  \\[2mm]
   & \le &
   \left( 1 + {1 \over p}\right)^{\! k} \,-\, 1  \;.
      \slabel{eq.prop.h_m(X,p).Delta.b}
\end{subeqnarray}
\end{corollary}

\noindent
Here the crude upper bound \reff{eq.prop.h_m(X,p).Delta.b}
is obtained from \reff{eq.prop.h_m(X,p).Delta.a}
by replacing $(k-rp+p)^{-1}$ by $p^{-1}$.
The true large-$k$ asymptotic behavior of \reff{eq.prop.h_m(X,p).Delta.a}
is $(1+1/p)^k k^{-1} p(p+1) [1 + O(1/k)]$.

Taking $p=1$ in \reff{eq.prop.h_m(X,p).Delta.a} gives:

\begin{corollary}
  \label{prop.h_m(X).Delta}
Let $X\subseteq V(G)$ where $|X|=k\geq 1$. Then
\be
   \sum\limits_{m=0}^\infty \Delta^{-m} h_m(X)
   \;\le\;
   {2 (2^k-1) \over k+1}  \;.
\ee
\end{corollary}

Using a similar proof technique to that of Proposition
\ref{prop.h_m(X,p,r).Delta}, but using
Proposition~\ref{prop.forests.Lambda}
instead of Proposition~\ref{prop.forests.Delta}, we may deduce

\begin{proposition}
  \label{prop.h_m(X,p,r).Lambda}
Let $X\subseteq V(G)$ where $|X|=k\geq rp$. Then
\be
   \sum\limits_{m=0}^\infty (m+r)^{-(k-1)}\Lambda^{-m} h_m(X,p,r)
   \;\le\;    rp^{-(r-1)}(k-rp+p)^{-1} {{k}\choose{r}}\;.
\ee
\end{proposition}

\section{Counting Connected Subgraphs (and Related Objects)} \label{sec.conn}



For $X \subseteq V(G)$, let $\scrc_m(X)$ be the set of all
$m$-edge subgraphs $H$ in $G$ (connected or not) such that
\begin{quote}
\begin{itemize}
   \item[(C1)] $X \subseteq V(H)$, and
   \item[(C2)] each component of $H$ contains at least one vertex in $X$.
\end{itemize}
\end{quote}
Note in particular the following special cases:
\begin{itemize}
   \item  For any $X$, $\scrc_0(X)$ has as its single element
       the edgeless graph with vertex set $X$.
   \item  $\scrc_0(\emptyset)$ has as its single element the
   empty graph, and $\scrc_m(\emptyset) = \emptyset$ for $m \ge 1$.
   \item  $\scrc_m(\{x\})$ consists of the connected $m$-edge subgraphs
       that contain $x$.
\end{itemize}
Recall the definition $w(H)=\prod_{e\in E(H)} w_e$.
Define the weighted counts
\be
   c_m(X)  \;=\;   \sum\limits_{H \in \scrc_m(X)}  w(H)  \;.
\ee
We then have the following bound in terms of maximum weighted degree:

\begin{proposition}
  \label{prop.conn}
Whenever $|X| = k\geq 0$ we have
\be
   c_m(X)  \;\le\; C(m,k) \, \Delta^m
\ee
where
\be
   C(m,k)   \;=\; \cases{ k (m+k)^{m-1} / m!     & for $k \neq 0$ \cr
              \noalign{\vskip 2mm}
              \delta_{m0}            & for $k=0$ \cr
            }
\ee
\end{proposition}

\noindent
The $k=1$ case of Proposition~\ref{prop.conn}
was proven previously \cite[Proposition~4.5]{Sokal_chromatic_bounds}.

Before beginning the proof of Proposition~\ref{prop.conn},
let us note some facts about the numbers $C(m,k)$:
\begin{enumerate}
\item  For each integer $m \ge 0$, $C(m,k)$ is a polynomial of degree $m$
in $k$.

\item  For each integer $m \ge 0$,
$C(m,k)$ is an increasing function of $k$ for $k \ge 0$.

\item  Generating function:  If ${\sf C}(z)$ solves the equation
\be
   {\sf C}(z)   \;=\;  e^{z {\sf C}(z)}
   \;,
 \label{gen_fn_eqn}
\ee
then
\be
   {\sf C}(z)^k \;=\;  \sum\limits_{m=0}^\infty C(m,k) \, z^m
 \label{gen_fn}
\ee
for all $k$ (integer or not);
this follows from the Lagrange inversion formula.\footnote{
   {\sc Proof.}  We use the Lagrange inversion formula in the form
   \cite[Theorem~5.4.2]{Stanley_99}:
   \begin{quote}
      Let $G(z) = \sum_{n=0}^\infty a_n z^n$ be a formal power series
      with $a_0 \neq 0$.  Then there is a unique formal power series
      $f(z) = \sum_{n=1}^\infty b_n z^n$ satisfying
      $f(z) = z \, G(f(z))$,
      and it satisfies
      $$  [z^n] f(z)^k  \;=\;  {k \over n} \, [z^{n-k}] G(z)^n  \;, $$
      where $[z^m] F(z)$ denotes the coefficient of $z^m$ in the
      formal power series $F(z)$.
   \end{quote}
   Now, equation~\reff{gen_fn_eqn} multiplied by $z$ has the given form
   with $f(z) = z {\sf C}(z)$ and $G = \exp$.  It follows that
   $$ [z^{n-k}] {\sf C}(z)^k  \;=\;
      {k \over n} \, [z^{n-k}] e^{nz}  \;=\;
      {k \over n} \, {n^{n-k} \over (n-k)!}   \;,$$
   which upon setting $n=m+k$ gives
   $ [z^m] {\sf C}(z)^k  =  k (m+k)^{m-1}/m! $.
}
Moreover, the series \reff{gen_fn} is absolutely convergent for $|z| \le 1/e$
and satisfies ${\sf C}(1/e) = e$.

\item For integer $k \ge 1$,
\be
   C(m,k)  \;=\;  \!\!\!\sum\limits_{\begin{scarray}
                    m_1, \ldots, m_k \ge 0 \\
                    m_1 + \cdots + m_k = m
                     \end{scarray}}
          \prod\limits_{i=1}^k C(m_i,1)
   \;.
     \label{cmk.identity1}
\ee
This is an immediate consequence of \reff{gen_fn}.

\item For all $k$ and $z$ (integer or not),
\be
   C(m,k)  \;=\;  \sum\limits_{f=0}^m {z^f \over f!} \, C(m-f,k-z+f)
   \;.
     \label{cmk.identity2}
\ee
This can easily be verified by direct calculation using the binomial formula:
\begin{eqnarray}
   & & \hspace*{-2.5cm}
   \sum\limits_{f=0}^m {z^f \over f!} \, {(k-z+f) (m+k-z)^{m-f-1} \over (m-f)!}
      \nonumber \\
   & = &
   {1 \over m!} \sum\limits_{f=0}^m {m \choose f} (k-z+f) z^f (m+k-z)^{m-f-1}
      \nonumber \\
   & = &
   {1 \over m!} \left( 1 \,-\, {d \over dk} \right)
        \sum\limits_{f=0}^m {m \choose f} z^f (m+k-z)^{m-f}
      \nonumber \\
   & = &
   {1 \over m!} \left( 1 \,-\, {d \over dk} \right) (m+k)^m
      \nonumber \\
   & = &
   {1 \over m!} \left[ (m+k)^m \,-\, m (m+k)^{m-1} \right]
      \nonumber \\
   & = &
   {k (m+k)^{m-1} \over m!}   \;.
\end{eqnarray}

\item For each fixed $k>0$, we have
\be
   C(m,k)   \;=\;  e^m m^{-3/2} {k e^k \over \sqrt{2\pi}} \, [1 + O(1/m)]
 \label{Cmk_asymptotics}
\ee
as $m \to\infty$.
This is an immediate consequence of Stirling's formula.
\end{enumerate}

\bigskip

In proving Proposition~\ref{prop.conn}, we will use the following two facts
about the weighted counts $c_m(X)$ for
$X=\{x_1,x_2,\ldots,x_k\}\subseteq V(G)$
where $x_1,x_2,\ldots,x_k$ are all distinct:

\bigskip
\noindent
{\bf Fact 1.} $c_m(X) \le \!\!\!
           \sum\limits_{\begin{scarray}
                  m_1, \ldots, m_k \ge 0 \\
                  m_1 + \cdots + m_k = m
                \end{scarray}}
           \! \prod\limits_{i=1}^k c_{m_i}(x_i) \;.$

\proof
Trivial when $k=1$.  For $k \ge 2$, construct a weight-preserving
bijection ${\bf F}$ of $\scrc_m(X)$ onto a subset of
$\!\!\!\!\bigcup\limits_{\begin{scarray}
                 m_1, \ldots, m_k \ge 0 \\
                 m_1 + \cdots + m_k = m
             \end{scarray}}
     \!\!\prod\limits_{i=1}^k \scrc_{m_i}(x_i)$
as follows:
Given $H \in \scrc_m(X)$, define ${\bf F}(H)_i$
to be the component of $H$ containing $x_i$ if this component
contains no vertex $x_{i'}$ with $i' < i$, and the  graph with
vertex set $\{x_i\}$ and no edges, otherwise.
\qed

\bigskip
\noindent
{\bf Fact 2.} $c_m(X) \le \sum\limits_{F \subseteq C(x_1)}
   w(F) \, c_{m-|F|}((X-x_1)\cup Y^F)$
where $C(x_1)$ is shorthand for the cocycle $E(\{x_1\}, \{x_1\}^c)$,
and $Y^F$ denotes the set of endpoints other than $x_1$ of the edges in $F$.
[Of course, the sum can be limited to sets $F$ having $|F| \le m$,
 since $c_m(X) = 0$ for $m < 0$.]

\proof
We classify the subgraphs $H \in \scrc_m(X)$
according to $F \equiv E(H)\cap C(x_1)$.
So, for each $F \subseteq C(x_1)$,
let ${\cal C}_m(X;F)$ be the set of all $H \in {\cal C}_m(X)$
that have $E(H) \cap C(x_1) = F$.
For $H\in {\cal C}_m(X;F)$, we let $H'=H \setminus x_1$.
Then $|E(H')|=m-|F|$
and $H'\in {\cal C}_{m-|F|}((X-x_1)\cup Y^F)$.
This gives an injective map from ${\cal C}_m(X;F)$ into
${\cal C}_{m-|F|}((X-x_1)\cup Y^F)$. Fact 2 now follows since
$w(H)=w(F)w(H')$.
\qed

We now give two alternative proofs of Proposition~\ref{prop.conn}.
The first proof uses Fact 2 for $k=1$ only, together with Fact 1;
it leads to a {\em nonlinear}\/ recursion whose solution is $C(m,1)$,
namely \reff{cmk.identity2} with $k=z=1$ combined with \reff{cmk.identity1}.
The second proof uses Fact 2 for all $k$, but does not use Fact 1;
it leads to a {\em linear}\/ recursion whose solution is $C(m,k)$,
namely \reff{cmk.identity2} with $z=1$.

We will need the following elementary result:

\begin{lemma}
  \label{subsetweight}
Let $S$ be a set in which each element $e\in S$
is given a nonnegative real weight $w_e$.
Then, for each integer $f \ge 0$, we have
\be
   \sum\limits_{\begin{scarray}
           F \subseteq S \\
           |F| = f
        \end{scarray}}
   \prod\limits_{e \in F} w_e
   \;\le\;
   {1 \over f!} \, \left( \sum\limits_{e \in S} w_e \right) ^{\! f}
   \;.
\ee
\end{lemma}

\par\medskip\noindent{\sc First proof of Proposition~\ref{prop.conn}.\ }
Consider first the case $k=1$, hence $X = \{x\}$.
We use induction on $m$. The proposition holds trivially when $m=0$,
so let us assume that $m\geq 1$.
We apply Fact 2 with $k=1$, and observe that the term $F=\emptyset$
contributes zero when $m \ge 1$  [since $c_{m-|F|}(Y^F)=c_{m}(\emptyset)=0$];
we therefore have
\be
   c_m(x) \;\le\;
          \!\!\!\!\sum\limits_{\begin{scarray}
                     \emptyset\neq F \subseteq C(x) \\
                     |F|\leq m
                       \end{scarray}}
          \! w(F) \; c_{m-|F|}(Y^F)
   \;.
\ee
Applying Fact 1 to bound $c_{m-|F|}(Y^F)$, we obtain
\be
   c_m(x) \;\le\;
   \!\!\!\!\sum\limits_{\begin{scarray}
              \emptyset\neq F \subseteq C(x) \\
              |F|\leq m
            \end{scarray}}
          w(F)
          \!\!\!\!\sum\limits_{\begin{scarray}
                     m_1, \ldots, m_{|Y^F|} \ge 0 \\
                     m_1 + \cdots + m_{|Y^F|} = m-|F|
                       \end{scarray}}
          \!\!\prod\limits_{i=1}^{|Y^F|} c_{m_i}(v_i)
\ee
where $Y^F = \{ v_1, \ldots, v_{|Y^F|} \}$.
Note that we have $m_i < m$ for all $i$ since $|F|\geq 1$.
We can therefore apply the inductive hypothesis
$c_{m_i}(v_i) \le C(m_i,1) \, \Delta^{m_i}$ to obtain
\begin{eqnarray}
   c_m(x)  & \le &
   \!\!\!\!\sum\limits_{\begin{scarray}
              \emptyset\neq F \subseteq C(x) \\
              |F|\leq m
            \end{scarray}}
          w(F)
          \!\!\!\!\sum\limits_{\begin{scarray}
                     m_1, \ldots, m_{|Y^F|} \ge 0 \\
                     m_1 + \cdots + m_{|Y^F|} = m-|F|
                       \end{scarray}}
          \!\!\prod\limits_{i=1}^{|Y^F|} C(m_i,1) \, \Delta^{m_i}
    \nonumber \\[2mm]
   & = &
   \!\!\!\!\sum\limits_{\begin{scarray}
              \emptyset\neq F \subseteq C(x) \\
              |F|\leq m
            \end{scarray}}
          w(F) \; C(m-|F|, \, |Y^F|) \; \Delta^{m-|F|}
    \nonumber \\[2mm]
   & \le &
   \!\!\!\!\sum\limits_{\begin{scarray}
              \emptyset\neq F \subseteq C(x) \\
              |F|\leq m
            \end{scarray}}
          w(F) \; C(m-|F|, \, |F|) \; \Delta^{m-|F|}
\end{eqnarray}
where the second line used the identity \reff{cmk.identity1},
and the last step used $|Y^F|\leq |F|$
and the fact that $C(m,k)$ is an increasing function of $k$.\footnote{
   If $G$ is a simple graph, we have $|Y^F| = |F|$.
   But if $G$ has multiple edges, all we can say is that $|Y^F|\leq |F|$.
}
Using Lemma~\ref{subsetweight}, we have
\begin{eqnarray}
   c_m(x)  & \le &
      \sum\limits_{f=1}^m  {\Delta^f \over f!} \, C(m-f,f) \, \Delta^{m-f}
    \nonumber \\[2mm]
   & = &
      \sum\limits_{f=0}^m  {\Delta^f \over f!} \, C(m-f,f) \, \Delta^{m-f}
    \nonumber \\[2mm]
   & = &  C(m,1) \, \Delta^m   \;,
\end{eqnarray}
where the second line used $C(m,0) = 0$ for $m \ge 1$,
and the last line used identity \reff{cmk.identity2} with $k=z=1$.
This proves Proposition~\ref{prop.conn} for $k=1$.

The result for general $k$ now follows using Fact 1
and the identity \reff{cmk.identity1}.
\qed

\par\medskip\noindent{\sc Second proof of Proposition~\ref{prop.conn}.\ }
We use induction on $m+k$.
The proposition holds trivially when $m=0$, so we assume $m\geq 1$.
Apply Fact 2 for general $k$;
because the right-hand side involves quantities $c_{m-|F|}((X-x_1)\cup Y^F)$
with $m-|F|+|(X-x_1)\cup Y^F|\leq m+k-1$,
we can apply the inductive hypothesis
$c_{m-|F|}((X-x_1)\cup Y^F) \le C(m-|F|, |(X-x_1)\cup Y^F|) \, \Delta^{m-|F|}$
to conclude that
\begin{eqnarray}
   c_m(X)  & \le &
      \sum\limits_{F \subseteq C(x_1)}
      w(F) \; C(m-|F|, \, |(X-x_1)\cup Y^F|) \; \Delta^{m-|F|}
   \nonumber \\[2mm]
   & \le &
      \sum\limits_{F \subseteq C(x_1)}
      w(F) \; C(m-|F|, \, k-1+|F|) \; \Delta^{m-|F|}
   \nonumber \\[2mm]
   & \le &
      \sum\limits_{f=0}^m
      {\Delta^f \over f!} \; C(m-f, \, k-1+f) \; \Delta^{m-f}
   \nonumber \\[2mm]
   & = & C(m,k) \, \Delta^m   \;,
\end{eqnarray}
where the second line used $|(X-x_1)\cup Y^F| \le k-1+|Y^F| \le k-1+|F|$
and the fact that $C(m,k)$ is increasing with $k$,
the third line used Lemma~\ref{subsetweight},
and the last line used identity \reff{cmk.identity2} with $z=1$.
\qed

\bigskip

Let us now show that Proposition~\ref{prop.conn} is best possible:

\bexam
  \label{example6.1}
Let ${\bf T}_r$ be the infinite $r$-regular tree,
let $x_1,\ldots,x_k \in V({\bf T}_r)$ satisfy
${\rm dist}(x_i,x_j) > m$ for $i \neq j$,
and let $G$ be the subtree of ${\bf T}_r$ induced by the vertices
lying at a distance at most $M$ from the set $\{x_1,\ldots,x_k\}$,
where $M$ is any fixed number greater than or equal to $m$.
Let all edge weights $w_e$ equal $\Delta/r$.
Then by a slight generalization of the computation in
\cite[proof of Proposition~4.2]{Sokal_chromatic_bounds}, one shows that
\be
   c_m(x_1,\ldots,x_k)   \;=\;
      \left( {\Delta \over r} \right) ^{\! m} \,
      {kr \over m} \,
      {kr + (r-1)m - 1 \choose m-1}
   \;,
\ee
which in the limit $r \to\infty$ with $\Delta$ fixed
tends to $C(m,k) \, \Delta^m$.
\eexam


\bexam
  \label{example6.2}
Let $G$ be the complete graph $K_n$,
and let all edge weights $w_e$ equal $\Delta/(n-1)$.
Consider first the case $k=1$:
let us count the subset of $\scrc_m(\{x\})$
consisting of the $m$-edge {\em trees}\/ containing $x$.
The number of such trees is ${n-1 \choose m} (m+1)^{m-1}$,
where the first factor counts the number of ways we can choose
$m$ additional vertices
and the second factor counts the number of trees on $m+1$ labelled vertices.
We therefore have
\be
   c_m(x)  \;\ge\;
  {n-1 \choose m} (m+1)^{m-1} \left({\Delta \over n-1}\right) ^{\! m}
  \;,
\ee
which in the limit $n \to\infty$ with $m$ fixed
tends to $C(m,1) \, \Delta^m$.

Now consider general $k$:
let $x_1,\ldots,x_k$ be distinct vertices,
and let us count the subset of $\scrc_m(x_1,\ldots,x_k)$
consisting of the $m$-edge $k$-component {\em forests}\/
in which each component contains exactly one $x_i$.
By a similar counting argument we have
\begin{eqnarray}
   & & \!\!\! c_m(x_1,\ldots,x_k)  \;\ge\;
         \nonumber \\
   & & \qquad
     \left({\Delta \over n\!-\!1}\right) ^{\! m}
     \!\!\!\!\sum\limits_{\begin{scarray}
                             m_1, \ldots, m_k \ge 0 \\
                             m_1 + \cdots + m_k = m
                          \end{scarray}}
  \prod_{i=1}^k  {n\!-\!k\!-\!m_1\!-\!\ldots\!-\!m_{i-1}  \choose  m_i}
                 \,   (m_i+1)^{m_i-1}
  \;, \qquad
  \nonumber \\[-5mm]
\end{eqnarray}
which in the limit $n \to\infty$ with $m$ fixed
tends to $C(m,k) \, \Delta^m$ [using \reff{cmk.identity1}].
\eexam

\bigskip

Note the difference between the ``pointwise'' bounds of
Propositions~\ref{prop4.1} and \ref{prop.conn}
and the ``generating-function'' bounds of
Propositions~\ref{prop4.2}, \ref{prop4.3}, \ref{prop.forests.Delta}
 and \ref{prop.forests.Lambda}.
In the former cases, the bound can be saturated
simultaneously for all (or arbitrarily many) $m$
--- as shown by the preceding two examples ---
so nothing can be gained by summing over $m$.
In the latter cases, by contrast,
the pointwise bound can be saturated only for one $m$ at a time
(cf.\ Example~\ref{sec.walks}.\ref{example4.3} after Corollary~\ref{cor4.5}),
so the generating-function bound is stronger.
The fundamental difference between the two situations
is that the subgraphs considered in
Propositions~\ref{prop4.1} and \ref{prop.conn}
are ``tied down only at one end''
(and hence can grow freely in all directions),
while those in
Propositions~\ref{prop4.2}, \ref{prop4.3}, \ref{prop.forests.Delta}
 and \ref{prop.forests.Lambda}
are ``tied down at all the leaves''.

\section{Counting Blocks, Block Paths, Block Trees and Block Forests}
   \label{sec.block}

We now seek an analogue of Proposition~\ref{prop.conn}
using maxmaxflow in place of maximum degree.
Unfortunately, no such bound is possible for the families $\scrc_m(X)$:
a simple counterexample is provided by the stars $G=K_{1,r}$,
which have maximum degree $r$ but maxmaxflow 1 (independent of $r$);
letting $x$ be the central vertex of the star,
we have $|\scrc_1(\{x\})| = r$, which is unbounded as $r\to\infty$.
The same thing happens for wheels $G = K_1 + C_r$,
so it is no help to assume that $G$ is non-separable.

To get a bound in terms of maxmaxflow, we need to restrict attention
to a suitably chosen proper {\em subfamily}\/ of $\scrc_m(X)$:
roughly speaking, we need to count subgraphs $H$ that possess
suitable ``non-separability properties''.
After much experimentation,
we settled on the family $\scrbt_m(X)$ of ``block trees''
and the families $\scrbf_m(X,Y)$ and $\scrbf_m^*(X,Y)$ of ``block forests'',
to be defined in Section~\ref{sec.block.1}.
However, simpler but weaker results can be obtained using
a different (and slightly larger) family $\scrb_m(X)$ of block forests,
to be defined in Section~\ref{sec.block.2}.
These two subsections are independent of each other
and can be read in either order.

\subsection[The classes $\scrbt_m(X)$, $\scrbf_m(X,Y)$ and $\scrbf_m^*(X,Y)$]{
            The classes \mbox{\boldmath$\scrbt_m(X)$},
                        \mbox{\boldmath$\scrbf_m(X,Y)$} and
                        \mbox{\boldmath$\scrbf_m^*(X,Y)$}}
  \label{sec.block.1}

Let $H$ be a not-necessarily-connected graph.
We recall that a vertex of $H$ is called an {\em internal vertex}\/ of $H$
if it is not a cut vertex of $H$.
Now let $B$ be a block of $H$.
We denote by ${\rm Int}(B,H)$
the set of internal vertices of $H$ that belong to $B$.
We say that $B$ is an {\em isolated block}\/ of $H$
if it contains no cut vertices of $H$,
and an {\em end block}\/ of $H$ if it contains exactly one cut vertex of $H$.
Finally, if $x,y \in V(H)$ with $x \neq y$,
we say that $H$ is an $xy$-{\em block path}\/
if $H$ is connected and is either a block (hence an isolated block of itself)
or else has exactly two end blocks $B_1,B_2$
with $x\in {\rm Int}(B_1,H)$ and $y\in {\rm Int}(B_2,H)$.

For any nonempty $X \subseteq V(G)$, let $\scrbt_m(X)$ be the set of all
$m$-edge subgraphs $H$ in $G$ such that
\begin{quote}
\begin{itemize}
   \item[(BT1)] $X \subseteq V(H)$;
   \item[(BT2)] $H$ is connected;
   \item[(BT3)] each end block $B$ of $H$
      contains at least one vertex of $X$ as an internal vertex
      [that is, ${\rm Int}(B,H)\cap X\neq \emptyset$]; and
   \item[(BT4)] if $H$ is a block, then either $H$ is an isolated vertex
      [hence $V(H)=X=\{x\}$]
      or else $H$ contains at least two vertices of $X$.
\end{itemize}
\end{quote}
By analogy with the set $\scrt_m(X)$,
which consists of trees whose leaves are ``tied down'' on the set $X$,
we think of the elements of $\scrbt_m(X)$ as {\em block trees}\/
whose end blocks are ``tied down'' on the set $X$.
[In particular, we have $\scrt_m(X) \subseteq \scrbt_m(X)$.]
Note the following special cases:
\begin{itemize}
   \item  $\scrbt_0(\{x\})$ has as its single element the edgeless graph
       with vertex set $\{x\}$,
       and $\scrbt_m(\{x\}) = \emptyset$ for $m \ge 1$.
   \item  When $x \neq y$,
       $\scrbt_m(\{x,y\})$ is the set of $m$-edge $xy$-block paths.
   \item  $\scrbt_0(X) = \emptyset$ for $|X| \ge 2$.
\end{itemize}

For $X,Y \subseteq V(G)$ with $Y\neq \emptyset$,
let $\scrbf_m(X,Y)$ be the set of all
$m$-edge subgraphs $H$ in $G$ (connected or not) such that
\begin{quote}
\begin{itemize}
   \item[(BF1)] $X\cup Y \subseteq V(H)$;
   \item[(BF2)] each component of $H$ contains exactly one vertex of $Y$;
   \item[(BF3)] each end block $B$ of $H$
      contains at least one vertex of $X \cup Y$ as an internal vertex
      [that is, ${\rm Int}(B,H)\cap (X \cup Y) \neq \emptyset$]; and
   \item[(BF4)] each isolated block of $H$ is either an isolated vertex
      belonging to $Y$ or else contains at least two vertices of $X\cup Y$.
\end{itemize}
\end{quote}
[Note that $\scrbf_m(X,Y) = \scrbf_m(X \setminus Y,Y)$,
so we can assume without loss of generality, if desired,
that $X \cap Y = \emptyset$.]
By analogy with the set $\scrf_m(X,Y)$,
which consists of forests whose leaves are ``tied down'' on the set $X\cup Y$
and whose components are ``tied down'' on single elements of the set $Y$,
we think of the elements of $\scrbf_m(X,Y)$ as {\em block forests}\/
whose end blocks are ``tied down'' on the set $X\cup Y$
and whose components are ``tied down'' on single elements of the set $Y$.
[In particular, we have
$\scrf_m(X,Y) \subseteq \scrbf_m(X,Y) \subseteq \scrc_m(Y)$.]
Note the following special cases:
\begin{itemize}
   \item  If $m=0$,
      then $\scrbf_0(X,Y) = \emptyset$ whenever $X\not\subseteq Y$, and
      $\scrbf_0(X,Y)$ has as its single element the edgeless graph with
      vertex set $Y$ whenever $X\subseteq Y$.
   \item  If $X = \emptyset$ (or more generally if $X\subseteq Y$),
      then $\scrbf_0(X,Y)$ has as its single element the edgeless graph with
      vertex set $Y$, and $\scrbf_m(X,Y) = \emptyset$ for $m \ge 1$.
   \item  If $X = \{x\}$ with $x \notin Y$, then each $H \in \scrbf_m(\{x\},Y)$
      is the disjoint union of an $m$-edge $xy$-block path for some $y\in Y$
      (this component avoiding the set $Y \setminus y$)
      and the collection $Y \setminus y$ of isolated vertices.
   \item  For $Y = \{y\}$, we have $\scrbf_m(X,\{y\}) = \scrbt_m(X \cup \{y\})$.
\end{itemize}

Recall that $w(H)=\prod_{e\in E(H)} w_e$.
Define the weighted counts
\begin{eqnarray}
   bt_m(X)    & = &   \sum\limits_{H \in  \scrbt_m(X)}  w(H)    \\[2mm]
   bf_m(X,Y)  & = &   \sum\limits_{H \in  \scrbf_m(X,Y)}  w(H)
\end{eqnarray}
We will obtain the following bounds on $bf_m(X,Y)$ in terms of
$\Delta$ and $\Lambda$:

\begin{proposition}
  \label{prop.blockforests.Delta}
For all $X,Y\subseteq V(G)$ with $Y \neq \emptyset$, we have
\be
   \sum\limits_{m=0}^\infty \left({\Delta \over \ln 2}\right)^{\! -m}
 bf_m(X,Y)
   \;\le\; 1   \;.
\ee
\end{proposition}

\begin{proposition}
  \label{prop.blockforests.Lambda}
For all $X,Y\subseteq V(G)$ with $Y \neq \emptyset$
and all $\alpha \in (1,2]$, we have
\be
   \sum\limits_{m=0}^\infty
\left({\alpha\Lambda \over \ln \alpha}\right)^{\! -m} bf_m(X,Y)
   \;\le\; \alpha^{|Y|-1}   \;.
\ee
\end{proposition}

Using the identity $\scrbf_m(X,\{y\}) = \scrbt_m(X \cup \{y\})$,
we immediately obtain the following corollaries for block trees:

\begin{corollary}
  \label{cor.blocktrees.Delta}
For all nonempty $X \subseteq V(G)$, we have
\be
   \sum\limits_{m=0}^\infty \left({\Delta \over \ln 2}\right)^{\! -m} bt_m(X)
   \;\le\; 1   \;.
\ee
\end{corollary}

\begin{corollary}
  \label{cor.blocktrees.Lambda}
For all nonempty $X \subseteq V(G)$, we have
\be
   \sum\limits_{m=0}^\infty \left({2\Lambda \over \ln 2}\right)^{\! -m} bt_m(X)
   \;\le\; 1   \;.
\ee
\end{corollary}

\noindent
Corollary~\ref{cor.blocktrees.Delta} could be proved
directly using the same proof technique as for Proposition
\ref{prop.blockforests.Delta}.
It is curious to note, however, that we have been unable to find a direct
proof of Corollary~\ref{cor.blocktrees.Lambda};
our proof of Proposition~\ref{prop.blockforests.Lambda} employs an inner
induction on $|Y|$, and thus inevitably passes through disconnected graphs.

As an immediate consequence of Corollary~\ref{cor.blocktrees.Lambda},
we obtain:

\begin{corollary}
  \label{cor.block.1}
Fix a weighted graph $(G,{\bf w})$ and an edge $e \in E(G)$.
Let $b_m(e)$ be the sum of the weights of the
non-separable $m$-edge subgraphs of $G$ containing $e$.
Then $b_1(e) = w_e$ and
\be
   \sum\limits_{m=2}^\infty
   \left({2\Lambda(G-e,{\bf w}) \over \ln 2}\right)^{\! -(m-1)} b_m(e)
   \;\le\; w_e   \;.
\ee
\end{corollary}
\bigskip\noindent
{\sc Proof of Corollary~\ref{cor.block.1},
     assuming Corollary~\ref{cor.blocktrees.Lambda}.}
Let $e = xy$ and put $X=\{x,y\}$. Clearly $b_1(e) = w_e$.
Let $H$ be a subgraph of $G$ and let $m\geq 2$.
Then $H$ is a non-separable $m$-edge subgraph of $G$ containing $e$
if and only if $H-e \in \scrbt_{m-1}(X)$.
Thus Corollary~\ref{cor.blocktrees.Lambda} applied to $G-e$ gives
\begin{eqnarray}
   \sum\limits_{m=2}^\infty
   \left({2\Lambda(G-e,{\bf w}) \over \ln 2}\right)^{\! -(m-1)} b_m(e)
   & = &
   w_e \sum\limits_{m=2}^\infty
   \left({2\Lambda(G-e,{\bf w}) \over \ln 2}\right)^{\! -(m-1)} bt_{m-1}(X)
       \nonumber \\[1mm]
   & \le &  w_e \;.
 \label{eq.proof.cor.block.1}
\end{eqnarray}
\phantom{p}\hfill\qed

\bigskip

Our proof of Proposition~\ref{prop.blockforests.Delta}
combines ideas from the proofs of Propositions~\ref{prop.forests.Delta}
and \ref{prop.conn}:
we use an inner induction on $|X|$ as a substitute for
the ``point-to-set'' bound of Proposition~\ref{prop4.2},
and we use a ``cutting'' argument similar to that employed
in the proofs of Propositions~\ref{prop4.2} and \ref{prop.conn}
to handle the case $|X|=1$.
In the inductive step we shall use
the following analogue of Lemma~\ref{lemma.decompose}
to split a block forest with $k$ end blocks
into a block forest with $k-1$ end blocks and a block path:

\begin{lemma}
  \label{lemma.decompose.blockforest}
Let $G$ be a graph, let $X,Y\subseteq V(G)$ with $Y \neq \emptyset$,
let $x\in X \setminus Y$, and let $H \in \scrbf_m(X,Y)$.
Let $H_1$ be the convex hull of $(X \setminus x)\cup Y$ in $H$,
and let $H_2$ be the convex hull of $\{x\}\cup V(H_1)$ in $H \setminus E(H_1)$.
Then $H$ is the edge-disjoint union of $H_1$ and $H_2$;
and for some $i$ $(0 \le i \le m)$,
we have $H_1 \in \scrbf_i(X \setminus x,Y)$
and $H_2 \in \scrbf_{m-i}(\{x\},V(H_1))$.
Moreover, the map $H \mapsto (H_1,H_2)$ is an injection.
\end{lemma}

\noindent
Note the slight change of perspective from Lemma~\ref{lemma.decompose}:
here $H_2$ is not an $xy$-block path for some $y \in V(H_1)$,
but rather the union of such a block path with the collection
$V(H_1) \setminus y$ of isolated vertices.
In particular, we have $V(H_1) \subseteq V(H_2)$.
However, modulo this change, this decomposition reduces to that of
Lemma~\ref{lemma.decompose} in the special case where $H \in \scrf_m(X,Y)$.

\proofof{Proposition~\ref{prop.blockforests.Delta}} As noted above,
we can assume without loss of generality that $X \cap Y =
\emptyset$. Let $c=1/\ln 2$.  We shall show that \be \label{e7.4}
   \sum\limits_{m=0}^M  \left(c\Delta\right)^{-m} bf_m(X,Y)
   \;\le\; 1
\ee
for all $M\geq 0$,
by using an outer induction on $M$ and an inner induction on $|X|$.
The base case $M=0$ and $|X|$ arbitrary
holds by the first remark after the definition of $\scrbf_m(X,Y)$.
The case when $X=\emptyset$ and $M$ is arbitrary
holds by the second remark following the definition of $\scrbf_m(X,Y)$.
Hence we may suppose that $M\geq 1$ and $|X|\geq 1$.
Our inductive argument consists of two steps:
\begin{quote}
\begin{quote}
\begin{itemize}
   \item[Step 1.] Proof that if \reff{e7.4} holds for all $|X|$
      and all $M'$ with $0\leq M' < M$, then it also holds for
      $|X|=1$ and $M$.
      This step uses a ``cutting'' argument.
   \item[Step 2.] Proof that if \reff{e7.4} holds for all $|X'|$ with
      $1\leq |X'|<|X|$  and some given $M$,
      then it holds for  $|X|$ and the same $M$.
      This step uses Lemma~\ref{lemma.decompose.blockforest}.
\end{itemize}
\end{quote}
\end{quote}

\medskip

{\em Step 1.}\/ Suppose that \reff{e7.4} holds for all $|X|$ and all
$M'$ with $0\leq M' < M$. Now let $X=\{x\}$ (note that $x \notin Y$
by assumption) and consider a subgraph $H\in \scrbf_m(\{x\},Y)$ for
some $m$. Note that $H$ is the disjoint union of an $xy$-block path
$H'$ for some $y\in Y$ and the collection $Y \setminus y$ of
isolated vertices. In particular, $x$ is neither an isolated vertex
nor a cut vertex of $H$. Let $F$ be the set of edges of $H$ incident
with $x$, let $f=|F|$ ($\ge 1$ because $x$ is not isolated), and let
$U^F$ be the set of end-vertices of edges in $F$ distinct from $x$.
Then $H'$ remains connected under deletion of $x$ (because $x$ is
not a cut vertex).

We shall show that $H \setminus x \in \scrbf_{m-f}(U^F,Y)$.
We have $U^F\subseteq V(H' \setminus x)$, so either $H' \setminus x=y$ or
$H' \setminus x$ contains at least two vertices of $U^F\cup \{y\}$.
For each end block $B$ of $H'\setminus x$, we have
${\rm Int}(B,H'\setminus x)\cap (U_F\cup \{y\})\neq \emptyset$,
otherwise
$B$ would be an end block of $H$
with
${\rm Int}(B,H)\cap (\{x\}\cup Y)= \emptyset$.
Thus $H \setminus x \in \scrbf_{m-f}(U^F,Y)$.

Let $C$ be the set of edges of $G$ incident to $x$.
For each nonempty $F\subseteq C$, let $\scrbf_m(\{x\},Y;F)$ be the set of all
subgraphs $H \in \scrbf_m(\{x\},Y)$ such that $E(H)\cap C=F$.
Then the map $H \mapsto H \setminus x$ is an injection from
 $\scrbf_m(\{x\},Y;F)$ into $  \scrbf_{m-f}(U^F,Y)$,
and $w(H) = w(F) w(H \setminus x)$.
Thus
\be
 bf_m(\{x\},Y) \;\le\;
 \sum_{\emptyset\neq F\subseteq C} w(F) \, bf_{m-f}(U^F,Y)   \;.
\ee
Hence
\begin{eqnarray}
   \sum\limits_{m=0}^M \left(c\Delta\right)^{-m}
bf_m(\{x\},Y)      & = &
   \sum\limits_{m=0}^M
\left(c\Delta\right)^{-m}
      \sum\limits_{\emptyset\neq F\subseteq C} w(F) \, bf_{m-f}(U^F,Y)
\nonumber \\[2mm]
      & = &
\sum\limits_{\emptyset\neq F\subseteq C}w(F)
\left(c\Delta\right)^{-f}
\sum\limits_{j=0}^{M-f} \left(c\Delta\right)^{-j}bf_j(U^F,Y)
\nonumber \\[2mm]
      & \le &  \sum\limits_{\emptyset\neq F\subseteq C}w(F)
\left(c\Delta\right)^{-f}
\end{eqnarray}
by the outer inductive hypothesis on $M$. Using Lemma~\ref{subsetweight},
we deduce that
\begin{eqnarray}
   \sum\limits_{m=0}^M \left(c\Delta\right)^{-m}
bf_m(\{x\},Y)
 &\leq&
\sum\limits_{f=1}^\infty
\left(c\Delta\right)^{-f}
\sum\limits_{\begin{scarray}
           F \subseteq C \\
           |F| = f
        \end{scarray}}
w(F)
\nonumber \\[2mm]
&\leq&
\sum\limits_{f=1}^\infty
\left(c\Delta\right)^{-f}
{\Delta^f\over f!}
\nonumber \\
&=&
e^{1/c}-1 \;=\; 1
\nonumber
\end{eqnarray}
since $c=1/\ln 2$. This proves Step 1.

\medskip

{\em Step 2.}\/
Suppose that \reff{e7.4} holds for all $|X'|$ with $1\leq |X'|<|X|$
and some given $M$.
Then $|X|\geq 2$. Choose arbitrarily some $x \in X$
(recall again that $x \notin Y$).
Given $H\in \scrbf_m(X,Y)$,
by Lemma~\ref{lemma.decompose.blockforest}
we may decompose $H$ into $H_1\cup H_2$
where $H_1 \in \scrbf_i(X \setminus x,Y)$
and $H_2 \in  \scrbf_{m-i}(\{x\},V(H_1))$.
Applying the inductive hypothesis that \reff{e7.4} holds
both for $X' = X \setminus x$ and for $X' = \{x\}$
with the given $M$,
we may deduce that
\begin{eqnarray}
& & \hspace*{-5mm}
\sum\limits_{m=0}^M (c\Delta)^{-m} bf_m(X,Y)
 \nonumber \\[2mm]
& & =\;
   \sum\limits_{m=0}^M  (c\Delta)^{-m} \!
      \sum\limits_{H\in  \scrbf_m(X,Y)} w(H)
 \nonumber \\[2mm]
& & \leq\;
\sum\limits_{m=0}^M
\sum\limits_{i=0}^m (c\Delta)^{-i} \!
      \sum\limits_{H_1\in  \scrbf_i(X \setminus x,Y)} w(H_1) \,
  (c\Delta)^{-(m-i)}
      \sum\limits_{H_2\in  \scrbf_{m-i}(\{x\},V(H_1))} w(H_2)
 \nonumber \\[2mm]
& & \leq\;
\sum\limits_{i=0}^M  (c\Delta)^{-i}
      \sum\limits_{H_1\in  \scrbf_i(X \setminus x,Y)} w(H_1)
\sum\limits_{j=0}^{M-i}  (c\Delta)^{-j}
\sum\limits_{H_2\in  \scrbf_{j}(\{x\},V(H_1))}  w(H_2)
 \nonumber \\
& & \leq\;  1  \;.
\end{eqnarray}
This completes the proof of Step 2 and hence of the proposition.
\qed

\medskip\noindent
{\bf Remarks.}  1.  Step 2 works for any value of $c$.
The specific value $c=1/\ln 2$ enters only in Step 1.

2. Step 2 is ``almost'' unnecessary;
we can ``almost'' apply Step 1 for any $X$.
The trouble is that if $|X| > 1$ we might find that
$x$ is a cut vertex of $H$.
In this case $H \setminus x \notin \scrbf_{m-f}(U^F,Y)$
since it will have one or more components containing no vertices of $Y$.

\medskip

\bexam
  \label{example7.1}
Let $G$ be the graph $K_2^{(s)}$ consisting of a pair of vertices $x,y$
joined by $s$ parallel edges.
Set all edge weights $w_e = \Delta/s$.
Let us consider $bt_m(\{x,y\}) = bf_m(\{x\},\{y\})$.
The generating function is
\be
   \sum_{m=0}^\infty \zeta^m bf_m(\{x\},\{y\})
   \;=\;
   \left( 1 + {\Delta \over s} \zeta \right) ^{\! s} \,-\, 1
   \;,
\ee which is an increasing function of $s$ (at fixed $\Delta$ and
$\zeta$) and tends to $e^{\Delta\zeta} - 1$ as $s \to \infty$. It
follows that Proposition~\ref{prop.blockforests.Delta} and
Corollary~\ref{cor.blocktrees.Delta} are sharp in the sense that
$(\ln 2)/\Delta$ is the maximal value of $\zeta$ that allows an
upper bound of 1. \eexam

\bexam
  \label{example7.2}
Let $G$ be the disjoint union of $n$ copies of $K_2^{(s)}$;
let $X$ contain one vertex from each copy,
and let $Y$ be the remaining vertices.
Set all edge weights $w_e = \Delta/s$.  Then
\be
   \sum_{m=0}^\infty \zeta^m bf_m(X,Y)
   \;=\;
   \left[ \left( 1 + {\Delta \over s} \zeta \right) ^{\! s} \,-\, 1
   \right] ^{\! n}
   \;\stackrel{s\to\infty}{\longrightarrow}\;
   (e^{\Delta\zeta} - 1)^n
   \;.
\ee So $(\ln 2)/\Delta$ is the maximal value of $\zeta$ that allows
{\em any}\/ finite upper bound that is independent of $|X|$ and
$|Y|$. \eexam

\medskip

Might it be possible to bound $\sum_{m=0}^\infty \zeta^m bf_m(X,Y)$
for some $\zeta > (\ln 2)/\Delta$ if we allow the right-hand side to
depend on $|X|$ and $|Y|$? We doubt it;  but all we can say for
sure, at present, is that $\zeta$ cannot exceed $(2\ln 2)/\Delta$:

\bexam
  \label{example7.3}
Let $G$ be the graph $P_n^{(s)}$ obtained from the $n$-edge path ($n \ge 2$)
by replacing each edge by $s$ parallel edges.
Set all edge weights $w_e = \Delta/(2s)$.
Let $x,y \in V(P_n^{(s)})$ with ${\rm dist}(x,y) = \ell$.  Then
\be
   \sum_{m=0}^\infty \zeta^m bf_m(\{x\},\{y\})
   \;=\;
   \left[ \left( 1 + {\Delta \over 2s} \zeta \right) ^{\! s} \,-\, 1
   \right] ^{\! \ell}
   \;\stackrel{s\to\infty}{\longrightarrow}\;
   (e^{\Delta\zeta/2} - 1)^\ell
   \;.
 \label{eq.example7.3}
\ee
Since $\ell$ can be arbitrarily large, a universal upper bound on
$\sum_{m=0}^\infty \zeta^m bf_m(X,Y)$ is impossible for $\zeta >
(2\ln 2)/\Delta$, even when $|X|=|Y|=1$. \eexam

%
%

\medskip
\bigskip

We have been unable to obtain a bound for $bf_m(X,Y)$ in terms of $\Lambda$
(i.e., Proposition~\ref{prop.blockforests.Lambda} or something like it)
by extending the proof technique of
Proposition~\ref{prop.forests.Lambda} in a similar way to the above proof.
The problem is that a universal ``point-to-set'' bound of the form
$\sum_{m=0}^\infty  (c\Lambda)^{-m} bf_m(X,Y)\leq 1$
(with a right-hand side that is independent of $|Y|$)
is simply not valid for any constant $c$:
it suffices to consider $G=K_{1,r}$ with $r>c$.
If, on the other hand, we try to adapt the proof technique of
Proposition~\ref{prop.blockforests.Delta}
by using an inductive hypothesis of the form
$\sum_{m=0}^\infty (m+|Y|)^{-(|X|-1)} (c\Lambda)^{-m} bf_m(X,Y) \leq |Y|$
(similar to that of  Proposition~\ref{prop.forests.Lambda}),
then we are unable to carry through Step 1
because of the increase in the size of $X$ when $\{x\}$ is replaced by
$U^F$.
Instead, our proof of Proposition~\ref{prop.blockforests.Lambda} will
rely solely on a ``cutting'' argument
rather than using Lemma~\ref{lemma.decompose.blockforest}.
In order for the induction to go through,
we need to work with a slightly larger
family of graphs than $ \scrbf_m(X,Y)$.

For $X,Y \subseteq V(G)$ with $Y\neq\emptyset$,
let $\scrbf_m^*(X,Y)$ be the set of all
$m$-edge subgraphs $H$ in $G$ (connected or not)
such that
\begin{quote}
\begin{itemize}
   \item[(BF1)] $X\cup Y \subseteq V(H)$;
   \item[(BF2*)] each component of $H$ contains at least one vertex of $Y$;
   \item[(BF3)] each end block $B$ of $H$
      contains at least one element of $X \cup Y$ as an internal vertex
      [that is, ${\rm Int}(B,H)\cap (X \cup Y) \neq \emptyset$]; and
   \item[(BF4)] each isolated block of $H$ is either an isolated vertex
      belonging to $Y$ or else contains at least two vertices of $X\cup Y$.
\end{itemize}
\end{quote}
The only change from $\scrbf_m(X,Y)$ is, therefore,
that each component of $H$ must contain {\em at least one}\/ vertex of $Y$,
rather than exactly one.
We have $\scrbf_m(X,Y) \subseteq \scrbf_m^*(X,Y) \subseteq \scrc_m(Y)$.
Since $\scrbf_m^*(X,Y) = \scrbf_m^*(X \setminus Y,Y)$,
we can assume without loss of generality, if desired,
that $X \cap Y = \emptyset$.
Note the special cases:
\begin{itemize}
   \item If $m=0$, then
      $\scrbf_0^*(X,Y) = \emptyset$ whenever $X\not\subseteq Y$, and
      $\scrbf_0^*(X,Y)$ has as its single element the edgeless graph with
      vertex set $Y$ whenever $X\subseteq Y$.
   \item If $Y=\{y\}$, then
      $\scrbf_m^*(X,\{y\}) = \scrbf_m(X,\{y\}) = \scrbt_m(X \cup \{y\})$.
\end{itemize}
Define the weighted counts
\be
   bf^*_m(X,Y)  \;=\;   \sum\limits_{H \in  \scrbf_m^*(X,Y)}  w(H)  \;.
\ee
Since $ \scrbf_m(X,Y) \subseteq  \scrbf_m^*(X,Y)$,
Proposition~\ref{prop.blockforests.Lambda} will follow from
the stronger result:

\begin{proposition}
  \label{prop.blockforestsdash.Lambda}
For all $X,Y\subseteq V(G)$ with $Y \neq \emptyset$
and all $\alpha \in (1,2]$, we have
\be
   \sum\limits_{m=0}^\infty
\left({\alpha\Lambda \over \ln \alpha}\right)^{\! -m} bf^*_m(X,Y)
   \;\le\; \alpha^{|Y|-1}   \;.
\ee
\end{proposition}

\proof Let $c=\alpha/\ln \alpha$.  We shall show that \be
   \sum\limits_{m=0}^M
\left(c\Lambda\right)^{-m} bf^*_m(X,Y)
   \;\le\; \alpha^{|Y|-1}   \;
\ee
for all $M\geq 0$,
by using induction on $M + |X| + |Y|$.
As noted above, we may assume that $X\cap Y=\emptyset$. The case
$M=0$ with $X,Y$ arbitrary holds by the remark after the definition
of $\scrbf_m^*(X,Y)$. The proposition also holds when $|X\cup Y|=1$
and $M$ is arbitrary, since $ bf^*_m(\emptyset,\{y\})=\delta_{m0}$.
So assume $|X \cup Y| \ge 2$ and choose $y \in Y$. By Proposition
\ref{prop2.cuts}, there exists a cocycle $C=E_G(L,R)$ in $G$ and a
vertex $z\in X \cup Y \setminus y$ such that $z\in L$, $X\cup Y
\setminus z \subseteq R$, and $\sum_{e \in C} w_e \le \Lambda$.
Later we shall distinguish two cases, depending on whether $z$
happens to lie in $X$ or in $Y$.

\begin{figure}[pt]
\vspace*{15mm}
\hspace*{5mm}
\input{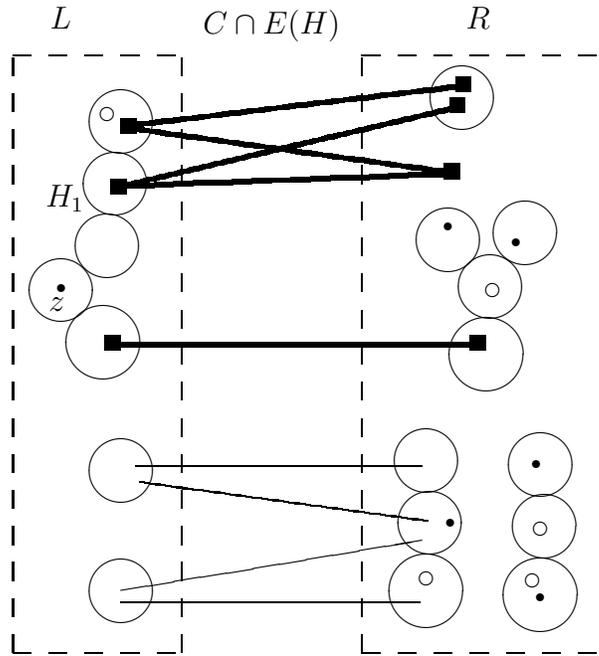}
\vspace{-3mm}
\caption{
   The graph $H \in \scrbf^*_m(X,Y)$.
   Blocks of $H_1$ and $H_2=H-H_1$ are either isolated vertices
   or else are indicated schematically by large open circles;
   the edges of $C\cap E(H)$ are shown explicitly, with those in
   $F$ drawn in bold.
   Vertices in $X$ are indicated by small solid circles,
   vertices in $Y$ are indicated by small open circles,
   and vertices in $L^F, R^F$ by small solid squares.
   Note that $H_2$ here contains one isolated vertex (which belongs to $R^F$)
   and that the four edges in $C\cap E(H)$ that do not belong $F$
   (which appear at the bottom of the figure) belong to $H_2$
   and do not generate vertices of $L^F$ or $R^F$.
}
 \label{fig.blockforestsdash.Lambda}
\end{figure}

For $H\in  \scrbf_m^*(X,Y)$, let $F:=F(H)$ be the set of edges of $H$
that occur as the first edge in $C$
on some path in $H$ from $z$ to $X\cup Y \setminus z$.
Let $H_1$ be the connected component of $H \setminus F$ containing $z$,
and let $H_2=H \setminus V(H_1)$.
Let $L^F$ (resp.\ $R^F$) be the set of vertices of $L$ (resp.\ $R$)
that are incident with $F$;
clearly $|L^F|,|R^F| \leq |F|$.
(This construction is illustrated in
 Figure~\ref{fig.blockforestsdash.Lambda}.)
For each nonempty $F\subseteq C$,
let $\scrbf^*_m(X,Y;F)$ be the set of all
subgraphs $H \in \scrbf^*_m(X,Y)$ such that $F(H) =F$.

Consider now the following two cases:

\bigskip

{\bf \boldmath Case 1: $z\in X$}.
Then $ \scrbf^*_m(X,Y;\emptyset)=\emptyset$.
We shall show that, for each nonempty $F\subseteq C$ with $|F|=f$
and each $H \in \scrbf_m(X,Y;F)$, we have
$H_1 \in  \scrbf_i^*(L^F,\{z\}) = \scrbt_i(L^F \cup \{z\})$ and
$H_2 \in  \scrbf_{m-f-i}^*(X \setminus z,Y\cup R^F)$
for some $i$ ($0\leq i\leq m-f$).
Note that by the above definitions $L^F\cup\{z\}\subseteq V(H_1)$
and $(X \setminus z)\cup (Y\cup R^F)\subseteq V(H_2)$.

We first consider $H_1$,
which is the component of $H\setminus F$ containing $z$.
Let $B_1$ be an end block  of $H_1$. Then $B_1$ is either
an end block of $H$, and hence $z\in {\rm Int}(B_1,H_1)$
or else $B_1$ is not an end block of $H$ and hence
${\rm Int}(B_1,H_1)\cap L^F\neq \emptyset$.
If $H_1$ is non-separable then either $V(H_1)=\{z\}$,
or else $H_1$ contains $z$ and at least one vertex of
$L^F\setminus z$ (otherwise $H_1$ would be an end block of $H$
with no internal vertex in $X \cup Y$).
Thus $H_1 \in  \scrbf_i^*(L^F,\{z\})$ for some $i$, $0\leq i\leq m-f$.

We next consider $H_2=H\setminus V(H_1)$.
Each component of $H_2$ is either a component of $H$, and hence
contains a vertex of $Y$, or is not a component of $H$, and hence
contains a vertex of $R^F$.
Let $B_2$ be an end block  of $H_2$. Then $B_2$ is either
an end block of
$H$, and hence satisfies
${\rm Int}(B_2,H_2)\cap [(X\setminus z)\cup Y]\neq \emptyset$,
or
is not an end block of $H$, and hence
satisfies ${\rm Int}(B_2,H_2)\cap R^F\neq \emptyset$.
Each isolated vertex of $H_2$ is either an isolated vertex of $H$,
and hence belongs to $Y$,
or is not an isolated vertex of $H$, and hence belongs to $R^F$.
Let $B_3$ be an isolated block of $H_2$ which is not a single vertex.
Then either: $B_3$ is an isolated block of $H$, and hence contains
at least two vertices of $(X\setminus z)\cup Y$; or $B_3$
is an end block of $H$, and hence contains
two distinct vertices,
one in $(X\setminus z)\cup Y$ and the other in $R^F$; or
$B_3$ is not a block of $H$, and hence contains
two distinct vertices
of $R^F$.
Thus $H_2 \in  \scrbf_{m-f-i}^*(X \setminus z,Y\cup R^F)$.

It follows that the  map $H \mapsto (H_1,H_2)$ is a
weight-preserving injection from \linebreak $\scrbf^*_m(X,Y;F)$ into
$\bigcup\limits_{i=0}^{m-f} \scrbf_i^*(L^F,\{z\}) \times
 \scrbf_{m-f-i}^*(X \setminus z,Y\cup R^F)$.
Therefore,
\be
 bf^*_m(X,Y)  \;\le\;  \sum_{\emptyset\neq F\subseteq C}w(F)
\sum_{i=0}^{m-f}bf_i^*(L^F,\{z\}) \; bf_{m-f-i}^*(X \setminus z,Y\cup R^F) \;.
\ee
It follows that
\begin{eqnarray}
  &&\!\!\!\! \sum\limits_{m=0}^M \left(c\Lambda\right)^{-m}
bf^*_m(X,Y)
\nonumber \\[2mm]
& \leq &
   \sum\limits_{m=0}^M
\left(c\Lambda\right)^{-m}
\sum_{\emptyset\neq F\subseteq C}w(F)
\sum_{i=0}^{m-f}bf_i^*(L^F,\{z\}) \; bf_{m-f-i}^*(X \setminus z,Y\cup R^F)
\nonumber \\[2mm]
      & = &
\sum\limits_{\emptyset\neq F\subseteq C}w(F)
\left(c\Lambda\right)^{-f}
\sum\limits_{i=0}^{M-f} \left(c\Lambda\right)^{-i}
bf_i^*(L^F,\{z\})
\sum\limits_{j=0}^{M-f-i} \left(c\Lambda\right)^{-j}
bf_j^*(X \setminus z,Y\cup R^F)
\nonumber \\[2mm]
      & \le &  \sum\limits_{\emptyset\neq F\subseteq C}w(F)
\left(c\Lambda\right)^{-f} \alpha^{|Y|-1+f}
\end{eqnarray}
by the inductive hypothesis on $M + |X| + |Y|$.
Using Lemma~\ref{subsetweight}, we deduce that
\begin{eqnarray}
   \sum\limits_{m=0}^M \left(c\Lambda\right)^{-m}
bf^*_m(X,Y)
 &\leq&
\sum\limits_{f=1}^\infty
\left(c\Lambda\right)^{-f} \alpha^{|Y|-1+f}
\sum\limits_{\begin{scarray}
           F \subseteq C \\
           |F| = f
        \end{scarray}}
w(F)
\nonumber \\[2mm]
&\leq&
\sum\limits_{f=1}^\infty
\left(c\Lambda\right)^{-f} \alpha^{|Y|-1+f} \,
{\Lambda^f\over f!}
\nonumber \\[2mm]
& = &  \alpha^{|Y|-1}(e^{\alpha/c}-1) \;\le\; \alpha^{|Y|-1}
\end{eqnarray}
since $c=\alpha/\ln \alpha$ and $\alpha \le 2$. Thus the proposition
holds when $z\in X$.

\bigskip

{\bf \boldmath Case 2: $z\in Y$}.
Then $y,z\in Y$ and hence $|Y|\geq 2$.
We can show in a similar way to Case 1 that,
for each $F\subseteq C$ with $|F|=f$
and each $H \in \scrbf^*_m(X,Y;F)$, we have
$ H_1 \in  \scrbf_i^*(L^F,\{z\})$ and
$ H_2 \in  \scrbf_{m-f-i}^*(X,(Y \setminus z)\cup R^F)$ for some $i$
($0\leq i\leq m-f$).
[Let us remark that if $F = \emptyset$,
 then $H_1$ is a component of $H$ consisting of just the
 isolated vertex $z$.]
The map $H \mapsto (H_1,H_2)$ is an injection from
$\scrbf^*_m(X,Y;F)$ into $\bigcup\limits_{i=0}^{m-f}
\scrbf_i^*(L^F,\{z\})\times
 \scrbf_{m-f-i}^*(X,(Y \setminus z)\cup R^F)$.
Thus
\be
 bf^*_m(X,Y)\leq \sum_{F\subseteq C}w(F)
\sum_{i=0}^{m-f} {bf}_i^*(L^F,\{z\}) \;
                 {bf}_{m-f-i}^*(X,(Y \setminus z)\cup R^F)   \;.
\ee
Hence
\begin{eqnarray}
&&\!\!\!\!   \sum\limits_{m=0}^M \left(c\Lambda\right)^{-m}
bf^*_m(X,Y)
\nonumber \\[2mm]
& \leq &
   \sum\limits_{m=0}^M
\left(c\Lambda\right)^{-m}
\sum_{F\subseteq C}w(F)
\sum_{i=0}^{m-f}{bf}_i^*(L^F,\{z\})\;
{bf}_{m-f-i}^*(X,(Y \setminus z)\cup R^F)
\nonumber \\[2mm]
      & = &
\sum\limits_{ F\subseteq C}w(F)
\left(c\Lambda\right)^{-f}
\sum\limits_{i=0}^{M-f} \left(c\Lambda\right)^{-i}
{bf}_i^*(L^F,\{z\})
\sum\limits_{j=0}^{M-f-i} \left(c\Lambda\right)^{-j}
{bf}_j^*( X,(Y \setminus z)\cup R^F)
\nonumber \\[2mm]
      & \le &  \sum\limits_{ F\subseteq C}w(F)
\left(c\Lambda\right)^{-f} \alpha^{|Y|-2+f}
\end{eqnarray}
by the inductive hypothesis on $M+|X|+|Y|$.
(Note that when $F=\emptyset$, we have $L^F = R^F = \emptyset$ and hence
$|(Y \setminus z)\cup R^F| = |Y|-1$.)
Using Lemma~\ref{subsetweight}, we deduce that
\begin{eqnarray}
   \sum\limits_{m=0}^M \left(c\Lambda\right)^{-m}
bf^*_m(X,Y)
 &\leq&
\sum\limits_{f=0}^\infty
\left(c\Lambda\right)^{-f} \alpha^{|Y|-2+f}
\sum\limits_{\begin{scarray}
           F \subseteq C \\
           |F| = f
        \end{scarray}}
w(F)
\nonumber \\[2mm]
&\leq&
\sum\limits_{f=0}^\infty
\left(c\Lambda\right)^{-f} \alpha^{|Y|-2+f} \,
{\Lambda^f\over f!}
\nonumber \\
& = &  \alpha^{|Y|-2} e^{\alpha/c}  \;=\;  \alpha^{|Y|-1}
\end{eqnarray}
since $c=\alpha/\ln \alpha$. Thus the proposition holds when $z\in
Y$. \qed

\medskip\noindent
{\bf Remark.}
The proof technique of Case 1 can be used
to prove the induction hypothesis
\be
   \sum\limits_{m=0}^M
   (c\Lambda)^{-m} bf^*_m(X,Y) \;\le\; \alpha^{|Y|-1}
 \label{xxx}
\ee
whenever $e^{\alpha/c} \le 2$.
(In particular, it can handle the apparently best-possible values
 $c=1/\ln 2$ and $\alpha=1$.)
Likewise, the proof technique of Case 2 can be used to prove the
induction hypothesis \reff{xxx} whenever $e^{\alpha/c} \le \alpha$.
The trouble is that we need the same hypothesis to work for both
cases, since we don't know {\em a priori}\/ whether $z$ will lie in
$X$ or in $Y$. Therefore, the best we can do --- at least with this
proof technique
--- seems to be to choose $\alpha \in (1,2]$ and then set
$c=\alpha/\ln \alpha$.

\bexam
  \label{exam.blockforest.lambda.1}
Let $T$ be any tree (e.g.\ a long path would do),
and let $G$ be the graph $T^{(s)}$ obtained from $T$
by replacing each edge by $s$ parallel edges.
Let all edge weights $w_e$ equal $\Lambda/s$,
so that the maxmaxflow is $\Lambda$.
Let $x,y \in V(T^{(s)})$ with ${\rm dist}(x,y) = \ell$.
Then $bf_m(\{x\},\{y\}) = bf^*_m(\{x\},\{y\})$, and we have
\be
   \sum\limits_{m=0}^\infty  \zeta^m bf_m(\{x\},\{y\})
   \;=\;
   \left[ \left( 1 + {\Lambda \over s} \zeta \right) ^{\! s} \,-\, 1
   \right] ^{\! \ell}
   \;\stackrel{s\to\infty}{\longrightarrow}\;
   (e^{\Lambda\zeta} - 1)^\ell
   \;.
\ee Since $\ell$ can be arbitrarily large, a universal upper bound
on $\sum_{m=0}^\infty \zeta^m bf_m(X,Y)$ is impossible for $\zeta >
(\ln 2)/\Lambda$, even when $|X|=|Y|=1$. \eexam

 \bexam
  \label{exam.blockforest.lambda.3}
Let $G$ be the union of $k$ disjoint copies of $K_{1,r}^{(s)}$,
with all edge weights $w_e = \Lambda/s$.
Let $X$ be the central vertices and $Y$ the remaining vertices,
so that $|X|=k$ and $|Y|=kr$.  Then
\be
   \sum\limits_{m=0}^\infty  \zeta^m bf_m(X,Y)
   \;=\;
   r^k \left[ \left( 1 + {\Lambda \over s} \zeta \right) ^{\! s} \,-\, 1 \right]
    ^{\! k}
   \;\stackrel{s\to\infty}{\longrightarrow}\;
   r^k (e^{\Lambda\zeta} - 1)^k
\ee
and
\be
   \sum\limits_{m=0}^\infty  \zeta^m bf^*_m(X,Y)
   \;=\;
   \left[ \left( 1 + {\Lambda \over s} \zeta \right) ^{\! rs}  \,-\, 1
   \right] ^{\! k}
   \;\stackrel{s\to\infty}{\longrightarrow}\;
   (e^{r\Lambda\zeta} - 1)^k
   \;.
\ee Thus, if any universal upper bound is possible for $\zeta=(\ln
2)/\Lambda$, the right-hand side has to be at least
$(|Y|/|X|)^{|X|}$ for $bf$ and $(2^{|Y|/|X|} - 1)^{|X|}$ for $bf^*$.
\eexam

\bigskip

These examples suggest the following conjectures:

\begin{conjecture}
   \label{conj.blockforests.Lambda}
For all $X,Y\subseteq V(G)$ with $Y \neq \emptyset$, we have
\be
   \sum_{m=0}^\infty
   \left({\Lambda \over \ln 2}\right)^{\! -m}  bf_m(X,Y)
   \;\le\;
   |Y|^{|X|}  \;.
\ee
\end{conjecture}

\begin{conjecture}
   \label{conj.blockforestsdash.Lambda}
For all $X,Y\subseteq V(G)$ with $Y \neq \emptyset$, we have
\be
   \sum_{m=0}^\infty
   \left({\Lambda \over \ln 2}\right)^{\! -m}  bf^*_m(X,Y)
   \;\le\;
   2^{|Y|} - 1  \;.
\ee
\end{conjecture}

\begin{conjecture}[a special case of Conjecture~\ref{conj.blockforests.Lambda}
                                     or \ref{conj.blockforestsdash.Lambda}]
   \label{conj.blocktrees.Lambda}
For all nonempty $X \subseteq V(G)$, we have
\be
   \sum_{m=0}^\infty
   \left({\Lambda \over \ln 2}\right)^{\! -m}  bt_m(X)
   \;\le\;
   1  \;.
\ee
\end{conjecture}

%

\subsection[The class $\scrb_m(X)$]{The class \mbox{\boldmath$\scrb_m(X)$}}
  \label{sec.block.2}



We conclude by discussing a larger class of block forests,
which is roughly the block analogue of the class $\scrh_m(X)$
considered in Section~\ref{sec.trees.2}.
For $X \subseteq V(G)$, let $\scrb_m(X)$ be the set of all
$m$-edge subgraphs $H$ in $G$ (connected or not)
such that
\begin{quote}
\begin{itemize}
   \item[(B1)] $X \subseteq V(H)$;
   \item[(B2)] each end block of $H$
      contains at least one element of $X$ as an internal vertex
      [that is, ${\rm Int}(B,H)\cap X\neq \emptyset$]; and
   \item[(B3)] each isolated block of $H$ is either an isolated vertex
      belonging to $X$ or else contains at least two vertices of $X$.
\end{itemize}
\end{quote}
[It follows from (B2) and (B3) that each component of $H$
is either an isolated vertex belonging to $X$
or else contains at least two vertices of $X$.
In particular, $\scrb_m(X) \subseteq \scrc_m(X)$.
Note also that
$\scrbf_m(X,Y) \subseteq \scrbf_m^*(X,Y) \subseteq \scrb_m(X \cup Y)$.]
Note the following special cases:
\begin{itemize}
   \item  For any $X$, $\scrb_0(X)$ has as its single element
       the edgeless graph with vertex set $X$.
   \item  $\scrb_m(\emptyset) = \emptyset$ for $m \ge 1$.
   \item  $\scrb_m(\{x\}) = \emptyset$ for $m \ge 1$.
   \item  $\scrb_m(\{x,y\})$ is the set of $m$-edge $xy$-block paths
       when $x \neq y$ and $m\geq 1$.
       \hfill\break
       [Hence $\scrb_m(\{x,y\}) = \scrbt_m(\{x,y\}) =
               \scrbf_m(\{x\},\{y\}) = \scrbf_m^*(\{x\},\{y\})$.]
\end{itemize}
Recall that $w(H)=\prod_{e\in E(H)} w_e$.
Define the weighted counts
\be
   b_m(X)  \;=\;   \sum\limits_{H \in \scrb_m(X)}  w(H)  \;.
\ee
We then have the following bound in terms of maxmaxflow:

\begin{proposition}
  \label{prop.blockforest}
Whenever $|X| = k \ge 1$ we have
\be
   b_m(X)  \;\le\; B(m,k) \, \Lambda^m
\ee
where
\be
   B(m,k)   \;=\;  2^m C(m, {k-1 \over 2})
        \;=\;  \cases{ (k-1) (2m+k-1)^{m-1} / m!     & for $k \neq 1$ \cr
              \noalign{\vskip 2mm}
              \delta_{m0}            & for $k=1$ \cr
            }
\ee
\end{proposition}

\smallskip

For the case of greatest interest,
namely $\scrb_m(\{x,y\}) = \scrbf_m^*(\{x\},\{y\})$,
the bound of Proposition~\ref{prop.blockforest} behaves
roughly like $(2e\Lambda)^m$,
which is much worse than the bound $(2\Lambda/\ln 2)^m$
of Propositions~\ref{prop.blockforests.Lambda} and
\ref{prop.blockforestsdash.Lambda}.
However, as we shall see, the proof of Proposition~\ref{prop.blockforest}
is quite a bit simpler than that of
Proposition~\ref{prop.blockforestsdash.Lambda}.

Before beginning the proof of Proposition~\ref{prop.blockforest},
let us note some facts about the numbers $B(m,k)$,
which follow easily from the facts about the $C(m,k)$
already discussed in Section~\ref{sec.conn}:
\begin{enumerate}
\item For each integer $m \ge 0$, $B(m,k)$ is an increasing function of $k \ge 1$.

\item  Generating function:  If ${\sf C}(z)$ solves the equation
\reff{gen_fn_eqn}, then
\be
   {\sf C}(2z)^{(k-1)/2}  \;=\;  \sum\limits_{m=0}^\infty B(m,k) \, z^m
 \label{gen_fn_B}
\ee
This is an immediate translation of \reff{gen_fn}.

\item  For all $k_1, k_2, m$ we have
\be
   \sum\limits_{i=0}^m B(i,k_1) \, B(m-i,k_2)   \;=\;
   B(m, k_1 + k_2 - 1)
   \;.
  \label{bmk.identity1}
\ee
This is an immediate consequence of \reff{gen_fn_B}.

\item For all $k$,
\be
   B(m,k)  \;=\;  \sum\limits_{f=0}^m {1 \over f!} \, B(m-f,k-1+2f)
   \;.
     \label{bmk.identity2}
\ee
This is an immediate translation of the identity \reff{cmk.identity2}
with $z=1/2$.

\item For each fixed $k>1$, we have
\be
   B(m,k)   \;=\;  (2e)^m m^{-3/2} {(k-1) e^{(k-1)/2} \over \sqrt{8\pi}}
           \, [1 + O(1/m)]
 \label{Bmk_asymptotics}
\ee
as $m \to\infty$.
This is an immediate translation of \reff{Cmk_asymptotics}.
\end{enumerate}
\bigskip

\par\medskip\noindent{\sc Proof of Proposition~\ref{prop.blockforest}.\ }
The Proposition holds trivially when $k = 1$,
so we may assume that $k \geq 2$.
Choose $x_1,x_2\in X$ and
let $C$ be a minimum cut in $G$ separating $x_1$ from $x_2$,
so that $\sum_{e \in C} w_e \le \Lambda$.
Let $G_1$ be the component of $G-C$ containing $x_1$,
and let $G_2$ be the union of the remaining components of $G-C$.
By construction, we have $x_1\in V(G_1)$ and $x_2\in V(G_2)$,
and each edge in $C$ joins a vertex of $G_1$ to a vertex of $G_2$.
Let $X_i=X\cap V(G_i)$ and $k_i = |X_i|$ for $i=1,2$.
Since $x_1 \in X_1$ and $x_2 \in X_2$,
we have $1 \le k_i \le k-1$ for $i=1,2$;
and of course $k_1 + k_2 = k$.

We shall classify the subgraphs $H$ of $G$ belonging to ${\cal B}_m(X)$
according to $F \equiv E(H)\cap C$.
So, for each $F \subseteq C$, let us set $f=|F|$;
define $Y_i^F$ ($i=1,2$) to be the set of vertices of $G_i$ that are
incident with $F$ (note that $|Y_i^F| \le f$);
and finally, let ${\cal B}_m(X;F)$ be the set of all $H \in {\cal B}_m(X)$
that have $E(H) \cap C = F$.

\begin{figure}[pt]
\input{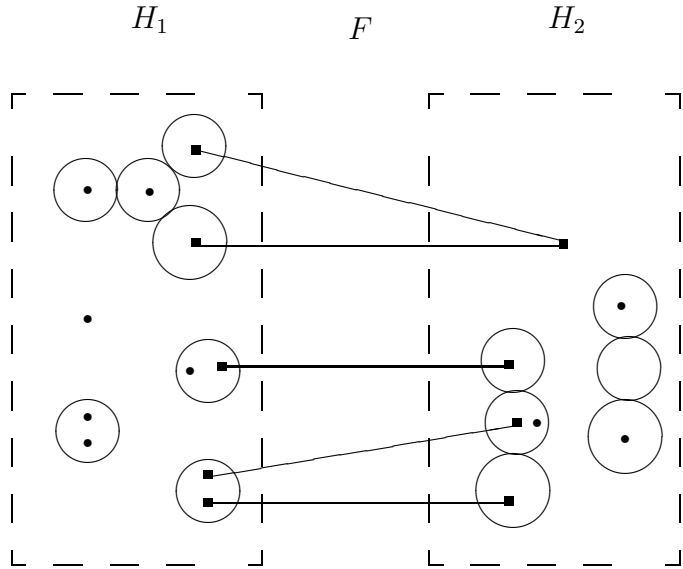}
\vspace{-5mm}
\caption{
   The graph $H \in {\cal B}_m(X)$.
   Blocks of $H_1$ and $H_2$ are indicated schematically by
   large open circles;
   the edges of $F$ are shown explicitly.
   Vertices in $X$ are indicated by small solid circles,
   and vertices in $Y_1^F, Y_2^F$ are indicated by small solid squares.
}
  \label{fig.blockforest}
\end{figure}

For $H\in {\cal B}_m(X;F)$, we let $H_i=H\cap G_i$ and $m_i = |E(H_i)|$;
note that $H = H_1 \cup H_2 \cup F$ and $m=m_1+m_2+f$.
We shall show that $H_i\in {\cal B}_{m_i}(X_i\cup Y_i^F)$ for $i=1,2$.
(The argument is illustrated in Figure~\ref{fig.blockforest}.)
Note first that since $X\subseteq V(H)$
we have $X_i\subseteq V(H_i)$.
Furthermore, each end block $B$ of $H_i$ is either an end block of $H$ or
else satisfies ${\rm Int}(B,H_i)\cap Y_i^F\neq \emptyset$.
Hence   ${\rm Int}(B,H_i)\cap (X_i\cup Y_i^F)\neq \emptyset$.
Each isolated vertex $v$ of $H_i$ is either an isolated vertex of $H$ or
else belongs to $Y_i^F$. Thus $v\in X_i\cup Y_i^F$.
Finally, each isolated block $B$ of $H_i$
that is not an isolated vertex of $H_i$
is either an isolated block of $H$, or else an end block of $H$ with its
cut vertex in $Y_i^F$, or else satisfies $|V(B)\cap Y_i^F|\geq 2$.
Thus $|V(B)\cap (X_i\cup Y_i^F)|\geq 2$.

It follows that we can construct a weight-preserving
[except for a factor $w(F)$] bijection
of $\scrb_m(X;F)$ onto a subset of
$\bigcup\limits_{i=0}^{m-f}
     \scrb_{i}(X_1\cup Y_1^F)\times \scrb_{m-f-i}(X_2\cup Y_2^F)$.
Thus
\be
b_m(X)  \;\leq\; \sum_{F\subseteq C}w(F)\sum_{i=0}^{m-f}
b_i(X_1\cup Y_1^F) \, b_{m-f-i}(X_2\cup Y_2^F)   \;.
\ee
Since $i \le m-f$, $|Y_1^F| \le f$ and $|X_1| < k$,
we have $i + |X_1 \cup Y_1^F| < m + k$;
and likewise we have $(m-f-i) + |X_2 \cup Y_2^F| < m + k$.
Therefore, we can use induction on $m+k$
(the Proposition being true for the initial case $m+k=1$)
to deduce that
\begin{eqnarray}
b_m(X)  &\leq&  \sum\limits_{F \subseteq C} w(F)
        \sum_{i=0}^{m-f} \Lambda^{m-f} \, B(i, |X_1 \cup Y_1^F|)
                           \, B(m-f-i, |X_2 \cup Y_2^F|)
   \nonumber \\[1mm]
    &\leq&  \sum\limits_{F \subseteq C} w(F)
        \sum_{i=0}^{m-f} \Lambda^{m-f} \, B(i, k_1+f)
                           \, B(m-f-i, k_2+f)
   \nonumber \\[1mm]
    &\leq&  \sum\limits_{f=0}^m {\Lambda^f \over f!}
        \sum_{i=0}^{m-f} \Lambda^{m-f} \, B(i, k_1+f)
                           \, B(m-f-i, k_2+f)
   \nonumber \\[1mm]
    & =  &  \Lambda^m B(m, k)   \;,
\end{eqnarray}
where the second line used the fact that $B(m,k)$ increases with $k$
for $k \ge 1$,
the third line used Lemma~\ref{subsetweight},
and the last line used
identities (\ref{bmk.identity1}) and (\ref{bmk.identity2})
and the fact that $k_1 + k_2 = k$.
\qed

\begin{corollary}
  \label{cor.block}
Fix a weighted graph $(G,{\bf w})$ and an edge $e \in E(G)$.
Then the sum of the  weights of the non-separable $m$-edge subgraphs of $G$
containing $e$ is at most $B(m-1,2) \, \Lambda(G-e,{\bf w})^{m-1}w_e$.
\end{corollary}

\proof
Let $e = xy$ and put $X=\{x,y\}$.
Let $H$ be a subgraph of $G$ and let $m\geq 1$.
Then $H$ is a non-separable $m$-edge subgraph of $G$ containing $e$
if and only if $H-e \in \scrb_{m-1}(X)$.
Thus Proposition~\ref{prop.blockforest} applied to $G-e$ gives
the claimed result.
\qed

\section*{Acknowledgments}

We wish to thank Jan van den Heuvel and Akira Saito
for valuable conversations.
We also thank Gordon Royle for drawing our attention to
references \cite{Bollobas_66,Mader_73},
and for many valuable conversations.
Finally, we thank an anonymous referee for several helpful comments,
including catching an error in our original versions of
Corollaries~\ref{cor.block.1} and \ref{cor.block}
(which were called Corollaries~7.5 and 7.12
 in the preprint version of this paper).

This research was supported in part by
U.S.\ National Science Foundation grants PHY--9900769, PHY--0099393
and PHY--0424082.
It was begun while one of the authors (A.D.S.)\ was a
Visiting Fellow at All Souls College, Oxford,
where he was hosted by the Department of Theoretical Physics
and supported in part by
U.K.\ Engineering and Physical Sciences Research Council grant GR/M 71626.

\end{document}